\documentclass[11pt]{article}
  \PassOptionsToPackage{dvipsnames}{xcolor}
  \usepackage[numbers,sort&compress]{natbib}
  \usepackage{amsthm}

\usepackage[utf8]{inputenc} %
\usepackage[T1]{fontenc}    %
\usepackage{microtype}      %

\usepackage{lipsum}         %
\usepackage{graphicx}       %
\usepackage{xcolor}         %
\usepackage{colortbl}       %
\usepackage{url}            %
\usepackage{nicefrac}       %
\usepackage{xspace}         %

\usepackage{amsmath}        %
\usepackage{amssymb}        %
\usepackage{bbm}            %
\usepackage{mathrsfs}       %
\usepackage{bm}             %
\usepackage{booktabs}       %
\usepackage{threeparttable} %
\usepackage{makecell}       %
\usepackage{subcaption}     %
\usepackage{multirow}       %
\usepackage{array}          %
\usepackage{tabularx}
\usepackage{arydshln}       %
\usepackage{setspace}       %

\makeatletter
\@ifpackageloaded{algorithm}{}{\usepackage{algorithm}}
\makeatother
\usepackage{algorithmic}    %
\usepackage{tikz}
\usetikzlibrary{arrows.meta,positioning,fit,shapes.misc}
\makeatletter
\@ifpackageloaded{hyperref}{}{\usepackage{hyperref}}
\makeatother
\usepackage{nameref}        %
\makeatletter
\@ifpackageloaded{cleveref}{}{\usepackage{cleveref}}
\makeatother

\hypersetup{
    colorlinks=true,
    linkcolor=blue,
    citecolor=blue,
    urlcolor=blue
}

\newcommand{\zero}{\bm{0}}
\newcommand{\RR}{\mathbb{R}}

\newcommand{\NN}{\mathbb{N}}

\newcommand{\prox}{\mathrm{prox}}

\newcommand{\calA}{\mathcal{A}}

\newcommand{\bB}{\bm{B}}

\newcommand{\bx}{\bm{x}}

\newcommand{\bu}{\bm{u}}
\newcommand{\bz}{\bm{z}}
\newcommand{\ba}{\bm{a}}
\newcommand{\by}{\bm{y}}

\newcommand{\bc}{\bm{c}}
\newcommand{\bU}{\bm{U}}

\newcommand{\blmbd}{\bm{\lambda}}

\newcommand{\bd}{\bm{d}}
\newcommand{\bq}{\bm{q}}
\newcommand{\bmu}{\bm{\mu}}
\DeclareMathOperator*{\argmin}{\mathrm{argmin}}

\DeclareMathOperator*{\argmax}{\mathrm{argmax}}

\newcommand{\calL}{\mathcal{L}}

\newcommand{\calC}{\mathcal{C}}

\newcommand{\calD}{\mathcal{D}}

\newcommand{\bv}{\bm{v}}
\newcommand{\bw}{\bm{w}}

\newcommand{\bV}{\bm{V}}

\newcommand{\bsigma}{\bm{\sigma}}

\makeatletter
\@ifundefined{theorem}{\newtheorem{theorem}{Theorem}[section]}{}

  \@ifundefined{lemma}{\newtheorem{lemma}{Lemma}[section]}{}
  \@ifundefined{corollary}{}{}
  \@ifundefined{proposition}{\newtheorem{proposition}{Proposition}[section]}{}
  \@ifundefined{definition}{\newtheorem{definition}{Definition}[section]}{}
  \@ifundefined{assumption}{\newtheorem{assumption}{Assumption}[section]}{}
  \@ifundefined{remark}{\newtheorem{remark}{Remark}[section]}{}
  \@ifundefined{example}{\newtheorem{example}{Example}[section]}{}

\@ifundefined{proof}{\newenvironment{proof}[1][Proof]{\par\noindent\textbf{#1.}\space}{\hfill$\square$\par}}{}
\@ifundefined{noproof}{}{}
\makeatother

\numberwithin{equation}{section}
\usepackage[normalem]{ulem}

\usepackage{amsfonts}
\ifpdf
  \DeclareGraphicsExtensions{.eps,.pdf,.png,.jpg}
\else
  \DeclareGraphicsExtensions{.eps}
\fi

\usepackage{enumitem}
\usepackage{placeins}       %
\setlist[enumerate]{leftmargin=.5in}
\setlist[itemize]{leftmargin=.5in}

\usepackage[most]{tcolorbox}
\newtcolorbox{stepbox}{enhanced, colback=gray!6, colframe=gray!55,
  boxrule=0.6pt, arc=2pt, left=6pt, right=6pt, top=4pt, bottom=4pt,
  before skip=8pt, after skip=6pt}

\newtcolorbox{questionbox}{enhanced, breakable, blanker,
  borderline west={2pt}{0pt}{gray!55},
  left skip=1em, right skip=1em, left=10pt, right=4pt,
  top=2pt, bottom=2pt, before skip=8pt, after skip=8pt}

  \usepackage[margin=1in]{geometry}
  \usepackage{mathpazo}
  \usepackage[scaled]{helvet}
  \usepackage{courier}
  \normalfont
  \allowdisplaybreaks[4]
  \hypersetup{
    colorlinks=true,
    citecolor=MidnightBlue,
    linkcolor=red!80!black,
    urlcolor=MidnightBlue,
    hypertexnames=false
  }

\makeatletter

  \@ifpackageloaded{lineno}{\AtBeginDocument{\nolinenumbers}}{}

\makeatother

\newcommand{\papertitle}{Dual-Accelerated Linearized Augmented Lagrangian
Method with a Constant Penalty for Nonconvex Composite Optimization}
\newcommand{\paperkeywords}{nonconvex composite optimization, augmented
Lagrangian method, single-loop method, constant penalty, dual regularization,
lower complexity bound}
\newcommand{\papermsc}{90C26, 90C30, 90C47, 90C60, 65K05, 49M37}
\newcommand{\paperabstract}{%
We study single-loop linearized augmented Lagrangian methods (ALMs) for
nonsmooth nonconvex composite optimization. Existing single-loop penalty methods
and linearized ALMs typically attain \(\mathcal O(\epsilon^{-3})\) complexity
for finding an \(\epsilon\)-KKT point by letting a penalty or regularization
parameter grow as \(\epsilon\) decreases. Whether a single-loop linearized ALM
can keep the penalty constant and still match or improve on this rate remains
unclear. By introducing a regularized augmented Lagrangian and constructing a
unified family of linearized regularized ALMs (LRALM) with moving centers, we
establish tight complexity bounds across constant center-update rules:
\(\Theta(\epsilon^{-4})\) when both the penalty and dual regularization are
accuracy-independent, and \(\Theta(\epsilon^{-3})\) when a constant penalty is
paired with \(\Theta(\epsilon^{-1})\) dual regularization. Even in the latter
regime, slow dual tracking prevents faster rates within this family. To improve
on this rate, we propose Fast LRALM, a single-loop method that exploits the
smooth strongly concave dual objective under a constant penalty via a Nesterov
extrapolation step on the multiplier update. Using one first-order oracle call
per iteration, Fast LRALM attains \(\mathcal O(\epsilon^{-5/2})\) iteration
complexity. To the best of our knowledge, it is the first ALM-based single-loop
method with a constant penalty to attain this rate for nonconvex composite
optimization, matching both the best-known complexity of multi-loop methods and
the optimal game-stationarity complexity for smooth nonconvex-concave minimax
problems.
}

  \title{\papertitle}
  \providecommand{\email}[1]{\href{mailto:#1}{\texttt{#1}}}
  \author{
  Luxuan Li\thanks{School of Mathematical Sciences, Beihang University,
  Beijing 100191, China; Faculty of Mathematics, University of Vienna,
   A-1090 Vienna, Austria (\email{lil00@univie.ac.at}).}
   \and
   Kangkang Deng\thanks{College of Science, National University of Defense
   Technology, Changsha 410073, China (\email{freedeng1208@gmail.com}).}
   \and
   Jiang Hu\thanks{Yau Mathematical Sciences Center, Tsinghua University,
   Beijing 100084, China (\email{jianghu@tsinghua.edu.cn}).}
   \and
   Chunfeng Cui\thanks{School of Mathematical Sciences and State Key
  Laboratory of Virtual Reality Technology and Systems, Beihang University,
   Beijing 100191, China (\email{chunfengcui@buaa.edu.cn}).}
  }
  \date{}

\begin{document}

\maketitle

  \begin{abstract}\paperabstract\end{abstract}

  \medskip
  \noindent\textbf{Mathematics Subject Classification.} \papermsc

\section{Introduction}
\label{sec:introduction}
In this paper, we study the nonsmooth nonconvex composite problem
\begin{align}
    \min_{\bx\in\RR^n}\quad & F(\bx) := f(\bx) + h(\calA(\bx)),
    \tag{\textbf{P}}
    \label{prob:main}
\end{align}
where \(f\colon\RR^n\to\RR\) is a possibly nonconvex function with a Lipschitz
continuous gradient, \(\calA\colon\RR^n\to\RR^m\) is a smooth and possibly nonlinear
mapping, and \(h\colon\RR^m\to\RR\) is a Lipschitz continuous convex function
whose proximal mapping is easy to evaluate; the complete assumptions are stated in
Section~\ref{sec:assumptions}. In this formulation, \(f\) is
usually a smooth loss or data-fidelity term, whereas \(h\) is a nonsmooth convex function
that imposes robustness, sparsity, or other structure on the transformed variable
\(\calA(\bx)\)~\cite{charisopoulos2021low}. Several models
of practical interest take this form.
\subsection{Motivating examples}

\textbf{Minimizing a pointwise maximum}~\cite{zhang2020smoothed}.
For smooth losses \(f_1,\dots,f_m\), such as logistic or robust classification
losses on \(m\) data subgroups, one solves
\begin{equation*}
  \min_{\bx\in\RR^n}\;\max_{1\le i\le m} f_i(\bx).
\end{equation*}
This is Problem~\eqref{prob:main} with \(h(\bu)=\max_{1\le i\le m}u_i\) and
\(\calA(\bx)=(f_1(\bx),\dots,f_m(\bx))\); a smooth regularizer, if present,
enters through \(f\).

\textbf{Nonconvex robust regression with structured regularization}~\cite{bohm2021variable}.
Given samples \(\{(\ba_j,b_j)\}_{j=1}^N\) and a linear analysis operator
\(\bB\in\RR^{m\times n}\), such as a finite-difference matrix for total variation
or fused lasso, one solves
\begin{equation*}
  \min_{\bx\in\RR^n}\;
  \frac{1}{2N}\sum_{j=1}^N\log\bigl(1+(\langle\ba_j,\bx\rangle-b_j)^2\bigr)
  +\|\bB\bx\|_1.
\end{equation*}
This is Problem~\eqref{prob:main} with the Cauchy loss \(f\),
\(\calA(\bx)=\bB\bx\), and \(h=\|\cdot\|_1\).

\textbf{Exact \(\ell_1\)-penalty reformulations}~\cite{burke1985descent,Drusvyatskiy2019}.
For a smooth constrained problem \(\min_{\bx\in\RR^n} f(\bx)\) subject to
\(\bc(\bx)\le\zero\), such as logistic or sigmoid risk constraints, an exact
penalty formulation solves
\begin{equation*}
  \min_{\bx\in\RR^n}\; f(\bx)+\nu\bigl\|[\bc(\bx)]_+\bigr\|_1,
\end{equation*}
with weight \(\nu>0\) and \([\bu]_+:=\max\{\bu,\zero\}\) applied componentwise.
This is Problem~\eqref{prob:main} with \(\calA=\bc\) and
\(h(\bu)=\nu\|[\bu]_+\|_1\).

\subsection{Related work}
A common approach to handling the nonsmooth composition \(h\circ\calA\) in
Problem~\eqref{prob:main} is to introduce an auxiliary variable \(\by=\calA(\bx)\)
and reformulate it as the nonlinear constrained problem
\begin{equation}\label{prob:split}
  \min_{\bx,\by}\; f(\bx)+h(\by)\quad\text{s.t.}\quad \calA(\bx)=\by.
\end{equation}
Most algorithms for Problems~\eqref{prob:main} and~\eqref{prob:split}, including
prox-linear~\cite{Drusvyatskiy2019,burke1985descent,lewis2016proximal},
quadratic penalty~\cite{kong2019complexity,lin2022complexity}, and inexact or
proximal augmented Lagrangian methods
(ALMs)~\cite{sahin2019inexact,li2021rate,xie2021proximal}, including
manifold  variants~\cite{deng2023manifold,deng2025manial,xu2025rial}, are
multi-loop methods that solve an inner subproblem approximately at each outer
iteration and require precise termination criteria. In large-scale settings,
single-loop designs become attractive. Existing single-loop penalty and
linearized ALM or alternating direction method of multipliers (ADMM)
methods~\cite{luSingleLoopGradientDescent2022,sun2024dual_descent_alm_admm,elbourkhissi2025perturbed},
however, typically attain an \(\mathcal O(\epsilon^{-3})\) complexity
for finding an \(\epsilon\)-KKT point by letting a penalty or regularization
parameter grow as the accuracy \(\epsilon\) decreases. As a result, the primal smoothness bound grows linearly
with the penalty parameter and restricts the primal stepsize accordingly,
which can slow convergence in practice. Keeping the penalty
constant, that is, independent of both the iteration counter and the
accuracy, avoids this growth and is a classical appeal of ALMs over penalty
methods~\cite{bertsekas1982,nocedal2006numerical}. In nonconvex optimization,
this appeal has been realized under linear constraints, where the smoothed
proximal ALM with a constant penalty attains
\(\mathcal O(\epsilon^{-2})\)~\cite{zhang2022global}. For general
nonlinear constraints, however, constant-penalty methods attain
\(\mathcal O(\epsilon^{-2})\) under uniform linear independence constraint
qualification~\cite{elbourkhissi2025linearized,liu2026fixed}
or by solving exact prox-linear subproblems under local
regularity~\cite{zhu2026nonconvex}. For orthogonality constraints, a
linearized ALM keeps the penalty constant but transfers
this parameter growth to a dual regularizer, so its complexity
remains \(\mathcal O(\epsilon^{-3})\)~\cite{zhu2026primal}. Beyond such special structures, Alacaoglu and
Wright~\cite{Alacaoglu24a} pose as an important open question whether an ALM with a
constant penalty and large dual stepsizes can attain
\(\mathcal O(\epsilon^{-3})\) for deterministic nonlinearly constrained
nonconvex problems, and a corresponding single-loop complexity theory remains
lacking for Problem~\eqref{prob:main} with a nonlinear
\(\calA\)~\cite{deng2026augmented}. For the minimax
reformulation of Problem~\eqref{prob:main}, multi-loop methods already attain
\(\mathcal O(\epsilon^{-5/2})\)~\cite{kong2021accelerated}, yet the best-known
single-loop ALM complexity remains \(\mathcal O(\epsilon^{-3})\). The
\(\Omega(\epsilon^{-2})\) lower bounds for smooth nonconvex
problems~\cite{carmon2020lower} and affinely constrained composite
problems~\cite{liu2026lower} do not determine whether this rate can be
improved, and existing lower bounds are established for individual single-loop
primal-dual algorithms under prescribed
stepsizes~\cite{mahdavinia2022tight,li2026smoothing}. This naturally raises the
following core question, which we aim to answer in this paper:
\begin{questionbox}
  \emph{Can a single-loop linearized ALM with a constant penalty improve on
  \(\mathcal O(\epsilon^{-3})\) for Problem~\eqref{prob:main}?}
\end{questionbox}
We answer it affirmatively in two steps, posed as
\hyperlink{Q1}{\textbf{(Q1)}} and \hyperlink{Q2}{\textbf{(Q2)}} below.

\subsection{Main contributions}
Writing \(h\) via its convex conjugate recasts Problem~\eqref{prob:main} as the
nonconvex-linear minimax formulation
\[
  \min_{\bx\in\RR^n}\max_{\bq\in\RR^m}\;
  \calL(\bx,\bq) := f(\bx)+\langle\bq,\calA(\bx)\rangle-h^*(\bq),
\]
whose game-stationary (GS) points coincide with the KKT points
of Problem~\eqref{prob:main}. On Riemannian manifolds, the single-loop RADA
framework~\cite{xu2026rada} exploits this nonconvex-linear structure
to reach \(\mathcal O(\epsilon^{-3})\) Riemannian game-stationarity
complexity through \(\epsilon\)-dependent dual regularization.
Alternatively, linearized ALM replaces the exact \(\bx\)-minimization of
classical ALM~\cite{hestenes1969multiplier,powell1969method,rockafellar1973multiplier}
by one gradient step on the reduced augmented
Lagrangian~\cite{dhingra2019proximal} and becomes gradient descent-ascent (GDA)
on a smooth nonconvex-concave (NC-C) minimax problem. For general smooth NC-C
minimax problems, primal smoothing~\cite{zhang2020smoothed} and dual
perturbation~\cite{lu2020hybrid,xu2023unified} each give
\(\mathcal O(\epsilon^{-4})\) in game stationarity, as does doubly smoothed
GDA~\cite{zheng2023universal}, while Perturbed Smoothed
GDA~\cite{li2026smoothing} reaches \(\mathcal O(\epsilon^{-3})\),
again through \(\epsilon\)-dependent dual regularization. In this paper,
we introduce the regularized augmented Lagrangian for the nonconvex composite
problem~\eqref{prob:main},
\[
  \mathcal P_s(\bx,\blmbd;\bz,\bmu)
  :=\calD_\rho(\bx,\blmbd)+\frac{1}{2\gamma}\|\bx-\bz\|^2
  -\frac{1}{2s}\|\blmbd-\bmu\|^2,
\]
where \(\calD_\rho\) is the reduced augmented Lagrangian with penalty parameter $\rho > 0$, \(\gamma,s>0\)
are regularization parameters, and \(\bz,\bmu\) are primal and dual centers;
see Section~\ref{subsec:regularized-al}.
This regularized formulation leads to our first question:
\begin{questionbox}
  \phantomsection\hypertarget{Q1}{\textbf{(Q1)}} \emph{With a constant
  penalty \(\rho\), can a linearized ALM attain
  \(\mathcal O(\epsilon^{-3})\) complexity for Problem~\eqref{prob:main}
  without growing the dual regularization parameter \(s\), and can growing
  \(s\) alone improve on \(\mathcal O(\epsilon^{-3})\)?}
\end{questionbox}

To address \hyperlink{Q1}{\textbf{(Q1)}}, we construct a unified family
\(\calC_\theta\) of linearized regularized ALMs (LRALM), presented in
Algorithm~\ref{alg:lralm}, that alternate GDA steps on \(\mathcal P_s\)
in \((\bx,\blmbd)\) with updates of the centers \((\bz,\bmu)\)
parameterized by constant \((\alpha,\beta)\in[0,1]^2\).
This family unifies linearized ALM and, when applied to \(\calD_\rho\), Smoothed
GDA~\cite{zhang2020smoothed}, Perturbed GDA~\cite{lu2020hybrid,xu2023unified},
doubly smoothed GDA~\cite{zheng2023universal}, and Perturbed Smoothed
GDA. Specifically, Theorems~\ref{thm:lower-bound}
and~\ref{thm:lower-bound-fixed-rho} together with
Propositions~\ref{prop:upper-beta1} and~\ref{prop:upper-beta0} provide a tight
worst-case complexity characterization
\[
    \underbrace{\Theta(\epsilon^{-4})}_{\text{small }\epsilon(\rho+s)}
    \qquad\text{and}\qquad
    \underbrace{\Theta(\epsilon^{-3})}_{\rho=\Theta(1),\; s=\Theta(\epsilon^{-1})}
\]
across all constant center movements \((\alpha,\beta)\in[0,1]^2\). In
particular, accuracy-independent parameters in \(\calC_\theta\) remain at
\(\Omega(\epsilon^{-4})\), indicating that \(\epsilon\)-dependent scaling
is needed within the family to surpass \(\epsilon^{-4}\). Keeping the penalty \(\rho\)
constant shifts this \(\epsilon^{-1}\) scaling to \(s\) while preserving a
\(\Theta(1)\) primal smoothness modulus and stepsize. These bounds also extend
earlier algorithm-specific minimax lower
bounds~\cite{mahdavinia2022tight,li2026smoothing} to the family
\(\calC_\theta\), showing that Smoothed GDA and Perturbed Smoothed GDA
attain the tight worst-case orders within \(\calC_\theta\) in their
respective regimes. Even with a constant penalty and growing \(s\), slow
dual tracking restricts the movement of the primal center independently of the
dual center, so faster rates require updates beyond \(\calC_\theta\).

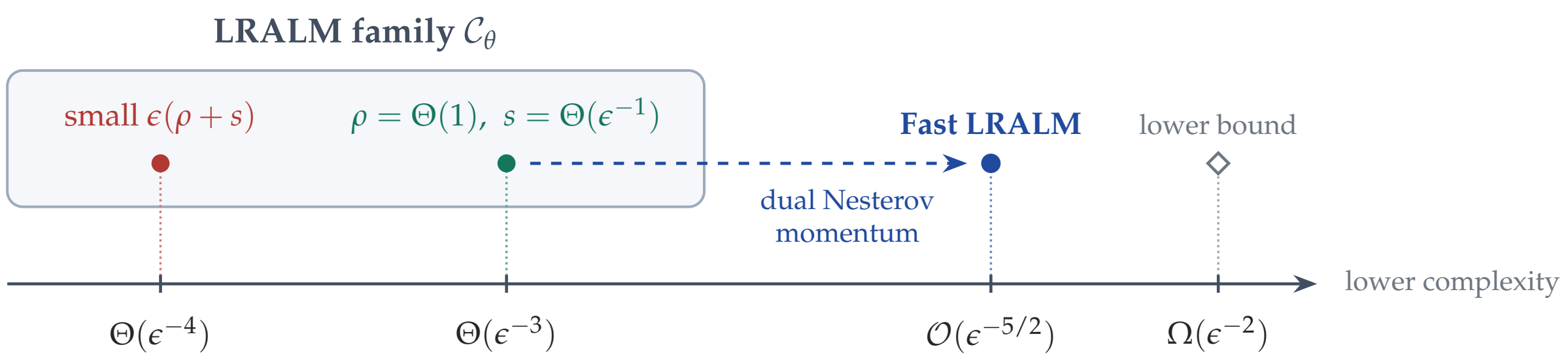
\begin{figure}[t]
\centering
\definecolor{cSlate}{RGB}{65,78,96}
\definecolor{cSlateBg}{RGB}{245,247,250}
\definecolor{cSlateBorder}{RGB}{156,170,188}
\definecolor{cRed}{RGB}{178,56,50}
\definecolor{cTeal}{RGB}{20,118,92}
\definecolor{cNavy}{RGB}{32,74,158}
\definecolor{cGray}{RGB}{110,118,128}
\resizebox{\ifdim\width>\textwidth\textwidth\else\width\fi}{!}{%
\begin{tikzpicture}[>=Stealth, font=\footnotesize]
  \def\xA{1.45}
  \def\xB{4.95}
  \def\xC{9.85}
  \def\xD{12.15}
  \def\yDot{1.22}

  \draw[cSlateBorder, line width=0.75pt, rounded corners=4.5pt, fill=cSlateBg]
    (-0.1, 0.78) rectangle (6.95, 2.16);
  \node[anchor=south, cSlate, font=\small\bfseries] at (3.425, 2.22)
    {LRALM family $\calC_\theta$};

  \draw[->, line width=0.85pt, cSlate] (-0.1, 0) -- (13.15, 0)
    node[right=4pt, font=\scriptsize, text=cGray] {lower complexity};

  \foreach \x/\col/\lab in {
      \xA/cRed/$\Theta(\epsilon^{-4})$,
      \xB/cTeal/$\Theta(\epsilon^{-3})$,
      \xC/cNavy/$\mathcal O(\epsilon^{-5/2})$,
      \xD/cGray/$\Omega(\epsilon^{-2})$
  }{
    \draw[densely dotted, line width=0.65pt, \col!65] (\x, 0.08) -- (\x, \yDot);
    \draw[line width=0.85pt, cSlate] (\x, -0.09) -- (\x, 0.09);
    \node[below=4pt, text=black!85] at (\x, -0.05) {\lab};
  }

  \draw[->, dashed, line width=0.9pt, cNavy]
    (\xB+0.24, \yDot) -- (\xC-0.24, \yDot);
  \node[below=3.5pt, font=\scriptsize, text=cNavy, align=center]
    at (8.40, \yDot) {dual Nesterov\\momentum};

  \fill[cRed] (\xA, \yDot) circle (2.6pt);
  \node[above=4.5pt, text=cRed] at (\xA, \yDot)
    {small $\epsilon(\rho+s)$};

  \fill[cTeal] (\xB, \yDot) circle (2.6pt);
  \node[above=4.5pt, text=cTeal] at (\xB, \yDot)
    {$\rho=\Theta(1)$,\; $s=\Theta(\epsilon^{-1})$};

  \fill[cNavy] (\xC, \yDot) circle (2.8pt);
  \node[above=4.5pt, text=cNavy, font=\footnotesize\bfseries] at (\xC, \yDot)
    {Fast LRALM};

  \draw[cGray, fill=white, line width=0.85pt]
    (\xD, \yDot) ++(0,3.0pt) -- ++(3.0pt,-3.0pt) -- ++(-3.0pt,-3.0pt) -- ++(-3.0pt,3.0pt) -- cycle;
  \node[above=4.5pt, text=cGray, font=\scriptsize] at (\xD, \yDot)
    {lower bound};
\end{tikzpicture}%
}
\caption{Summary of our iteration-complexity bounds for finding an
\(\epsilon\)-KKT point of Problem~\eqref{prob:main}, comparing the LRALM family
\(\calC_\theta\) and Fast LRALM against the \(\Omega(\epsilon^{-2})\) baseline
lower bound~\cite{carmon2020lower,liu2026lower}.}
\label{fig:complexity-axis}
\end{figure}

\begin{table}[t]
\centering
\begin{threeparttable}
\caption{Comparison of complexity bounds (iterations for single-loop methods;
first-order oracle calls for multi-loop methods) for finding an
\(\epsilon\)-GS point (smooth NC-C minimax) or an \(\epsilon\)-KKT point
(nonlinear constrained optimization).}
\label{tab:known-results}
\scriptsize
\setlength{\tabcolsep}{2.4pt}
\renewcommand{\arraystretch}{1.12}
\setlength{\aboverulesep}{0pt}
\setlength{\belowrulesep}{0pt}
\begin{tabular*}{\textwidth}{l @{\extracolsep{\fill}} lcccc}
\toprule
\multirow{2}{*}[-1.2ex]{Problem} & \multirow{2}{*}[-1.2ex]{Primal-dual algorithm}
 & \multirow{2}{*}[-1.2ex]{Single-loop}
 & \multicolumn{2}{c}{\rule{0pt}{2.2ex}Core parameters}
 & \multirow{2}{*}[-1.2ex]{Complexity} \\
\cmidrule(lr){4-5}
 & & & Penalty \(\rho\) & Dual reg.\ \(s\)
 & \\[2.5pt]
\midrule
\multirow{5}{*}{\makecell[l]{Smooth NC-C\\ minimax}}
 & Smoothed GDA~\cite{zhang2020smoothed,li2026smoothing}
 & \checkmark & n/a & --
 & \(\Theta(\epsilon^{-4})\) \\
 & Perturbed Smoothed GDA~\cite{li2026smoothing}
 & \checkmark & n/a & \(\Theta(\epsilon^{-1})\)
 & \(\Theta(\epsilon^{-3})\) \\
 & AIPP-S~\cite{kong2021accelerated}
 & -- & n/a & \(\Theta(\epsilon^{-1})\)
 & \(\mathcal O(\epsilon^{-5/2})\) \\
 & Tracked-FOAM~\cite{pan2026optimal}
 & -- & n/a & \(\Theta(\epsilon^{-1})\)
 & \(\Theta(\epsilon^{-5/2})\) \\
 & MC-EGM~\cite{zhang2026matching}
 & \checkmark & n/a & \(\Theta(\epsilon^{-1})\)
 & \(\Theta(\epsilon^{-5/2})\) \\
\midrule
\multirow{9}{*}{\makecell[l]{Nonlinear\\ constrained\\ optimization}}
 & iALM (Sahin et~al.)~\cite{sahin2019inexact}
 & -- & \(\Theta(\epsilon^{-1})\) & --
 & \(\tilde{\mathcal O}(\epsilon^{-4})\) \\
 & iALM (Li et~al.)~\cite{li2021rate}
 & -- & \(\Theta(\epsilon^{-1})\) & --
 & \(\tilde{\mathcal O}(\epsilon^{-3})\) \\
 & GDPA~\cite{luSingleLoopGradientDescent2022}
 & \checkmark & \(\Theta(\epsilon^{-1})\) & \(\Theta(\epsilon^{-1})\)
 & \(\mathcal O(\epsilon^{-3})\) \\
 & LIPAL~\cite{elbourkhissi2025perturbed}
 & \checkmark & \(\Theta(\epsilon^{-1})\) & \(\Theta(\epsilon^{-1})\)
 & \(\mathcal O(\epsilon^{-3})\) \\
 & LSALM~\cite{zhu2026primal}
 & \checkmark & \(\Theta(1)\) & \(\Theta(\epsilon^{-1})\)
 & \(\mathcal O(\epsilon^{-3})\) \\
 & NR-LALM~\cite{liu2026fixed}
 & \checkmark & \(\Theta(1)\) & --
 & \(\mathcal O(\epsilon^{-2})\) \\[2.5pt]
\cdashline{2-6}[2pt/2pt]
\noalign{\vskip 2pt}
 & \phantom{\textbf{LRALM} \(\boldsymbol{\calC_\theta}\)} (Thm.~\ref{thm:lower-bound}, Prop.~\ref{prop:upper-beta1})
 & \checkmark & \(\Theta(1)\) & \(\Theta(1)\)
 & \(\Theta(\epsilon^{-4})\) \\
 & \multirow{-2}{*}{\textbf{LRALM} \(\boldsymbol{\calC_\theta}\)} (Thm.~\ref{thm:lower-bound-fixed-rho}, Prop.~\ref{prop:upper-beta0})
 & \checkmark & \(\Theta(1)\) & \(\Theta(\epsilon^{-1})\)
 & \(\Theta(\epsilon^{-3})\) \\[1.5pt]
 & \cellcolor{gray!22}\rule{0pt}{2.3ex}\rule[-0.8ex]{0pt}{0pt}\textbf{Fast LRALM} (Alg.~\ref{alg:b2-prime})
 & \cellcolor{gray!22}\checkmark
 & \cellcolor{gray!22}\(\Theta(1)\) & \cellcolor{gray!22}\(\Theta(\epsilon^{-1})\)
 & \cellcolor{gray!22}\(\mathcal O(\epsilon^{-5/2})\) \\
\bottomrule
\end{tabular*}
\begin{tablenotes}[flushleft]
\scriptsize
\item[] \textbf{Notes.}\;
``--'' indicates that the property is not satisfied or the
parameter is not used, and ``n/a'' marks an item not applicable to smooth
minimax. AIPP-S applies to Problem~\eqref{prob:main} through the minimax
reformulation. For NR-LALM, each primal step solves a Gauss--Newton-type linear
system instead of taking an explicit gradient step. For LIPAL, the single-loop
entry refers to its prox-gradient variant. LSALM is restricted to orthogonality
constraints, where a local uniform constraint qualification holds automatically,
and requires a compact primal domain and dual
projection. For GDPA, LIPAL, and LSALM, \(s\) denotes the dual regularization
parameter obtained by writing the dual update in the form of \(\mathcal P_s\).
Lower bounds in \(\Theta(\cdot)\)
are algorithm-specific for Smoothed GDA and Perturbed Smoothed GDA, hold over
first-order method classes on smooth NC-C problems for Tracked-FOAM and MC-EGM,
and hold over constant center stepsizes in \(\calC_\theta\) for LRALM.
\end{tablenotes}
\end{threeparttable}
\end{table}

In smooth NC-C minimax optimization, faster rates have recently
been attained. The optimal deterministic first-order complexity
for \(\epsilon\)-GS points is \(\Theta(\epsilon^{-5/2})\), matched by the
multi-loop methods~\cite{kong2021accelerated,lin2020near,ostrovskii2021efficient,pan2026optimal} and by the
single-loop moving-center extragradient method with momentum (MC-EGM) of
Zhang et~al.~\cite{zhang2026matching}. While these methods provide strong
complexity guarantees for general smooth minimax problems,
Problem~\eqref{prob:main} carries additional structure through \(\calD_\rho\)
that a tailored ALM can exploit.
This motivates our second question:
\begin{questionbox}
  \phantomsection\hypertarget{Q2}{\textbf{(Q2)}} \emph{Beyond \(\calC_\theta\),
  how can a constant-penalty linearized ALM exploit this structure to overcome
  the \(\Theta(\epsilon^{-3})\) barrier and attain
  \(\mathcal O(\epsilon^{-5/2})\) for Problem~\eqref{prob:main}?}
\end{questionbox}

To answer \hyperlink{Q2}{\textbf{(Q2)}}, we develop Fast LRALM,
presented in Algorithm~\ref{alg:b2-prime}, which exploits the smooth
and strongly concave dual objective induced by \(\mathcal P_s\)
under a constant
penalty by adding a single Nesterov accelerated gradient (NAG) extrapolation
step~\cite{nesterov2018lectures} to the multiplier update. This update
accelerates dual tracking and relaxes the restriction on primal center movement
in \(\calC_\theta\), requiring one first-order oracle call per iteration
without inner loops. A tailored Lyapunov construction, combining NAG-type energy with
primal tracking and an \(\epsilon\)-independent multiplier bound from the
constant penalty, yields an
\(\mathcal O(\epsilon^{-5/2})\) iteration complexity for finding an
\(\epsilon\)-KKT point. To the best of our knowledge, Fast
LRALM is the first ALM-based single-loop method with a constant penalty to
attain \(\mathcal O(\epsilon^{-5/2})\) for Problem~\eqref{prob:main}, matching
both the best-known order of multi-loop methods~\cite{kong2021accelerated} and
the optimal order for game stationarity in smooth NC-C minimax
optimization~\cite{pan2026optimal,zhang2026matching}.
These results also provide partial progress toward the open question of
Alacaoglu and Wright~\cite{Alacaoglu24a} in this composite setting; see
Remark~\ref{rem:lralm-positioning}(ii) for a detailed comparison.
Figure~\ref{fig:complexity-axis} summarizes our results, and
Table~\ref{tab:known-results} compares our complexity bounds and parameter
scales with representative methods.

\paragraph{Organization.}
Section~\ref{sec:preliminaries} introduces the problem setup and preliminaries.
Section~\ref{sec:algorithm-framework} defines the LRALM family and establishes
its tight complexity bounds. Section~\ref{sec:acceleration} presents Fast LRALM
and its convergence analysis.
Section~\ref{sec:experiments} reports numerical experiments, and
Section~\ref{sec:conclusion} concludes the paper.

\paragraph{Notation.} Throughout the paper, \(\|\cdot\|\) denotes the spectral norm of a
matrix and the Euclidean norm of a vector, and \(\langle\cdot,\cdot\rangle\) stands for the Euclidean inner product. For
\(\calA=(A_1,\ldots,A_m)^\top\), the Jacobian \(J_{\calA}(\bx)\in\RR^{m\times n}\)
has \(i\)th row \(\nabla A_i(\bx)^\top\).
Standard asymptotic notation
\(\mathcal O(\cdot)\), \(\Omega(\cdot)\), and \(\Theta(\cdot)\) hides universal
positive constants independent of the accuracy \(\epsilon\) and the dimensions
\(n,m\).

\section{Preliminaries}
\label{sec:preliminaries}

We collect the problem setup and oracle model for
Problem~\eqref{prob:main}. We then introduce the regularized augmented
Lagrangian that underlies our algorithms. Following~\cite{rockafellar1970convex}, for a proper closed convex function
\(h\colon\RR^m\to\RR\cup\{+\infty\}\), its subdifferential, convex conjugate,
proximal mapping, and Moreau envelope are
\begin{gather*}
    \partial h(\by):=\{\bv:\ h(\bw)\ge h(\by)+\langle\bv,\bw-\by\rangle
    \ \text{for all }\bw\},
    \qquad
    h^*(\bq):=\sup_{\by\in\RR^m}\{\langle\bq,\by\rangle-h(\by)\},\\
    \prox_{t h}(\bu)
    :=\argmin_{\by\in\RR^m}
       \Bigl\{h(\by)+\frac{1}{2t}\|\by-\bu\|^2\Bigr\},
    \quad
    M_t^h(\bu)
    :=\min_{\by\in\RR^m}
       \Bigl\{h(\by)+\frac{1}{2t}\|\by-\bu\|^2\Bigr\},
    \quad t>0.
\end{gather*}

\subsection{Problem setup}
\label{sec:assumptions}
We impose the following standing assumptions on Problem~\eqref{prob:main}.
\begin{assumption}
\label{ass:smoothness}
There exist constants \(L_f,G_f,L_{\calA},G_{\calA}\ge0\) such that, for
all \(\bx,\bx'\in\RR^n\),
\[
\begin{aligned}
    \|\nabla f(\bx)-\nabla f(\bx')\|
    &\le L_f\|\bx-\bx'\|,
    &\|\nabla f(\bx)\|&\le G_f,\\
    \|J_{\calA}(\bx)-J_{\calA}(\bx')\|
    &\le L_{\calA}\|\bx-\bx'\|,
    &\|J_{\calA}(\bx)\|&\le G_{\calA}.
\end{aligned}
\]
\end{assumption}

\begin{assumption}
\label{ass:h}
The function \(h\colon\RR^m\to\RR\) is convex and \(L_h\)-Lipschitz, in the sense that
\[
    |h(\bu)-h(\bv)|\le L_h\|\bu-\bv\|,
    \qquad \bu,\bv\in\RR^m.
\]
\end{assumption}

We measure approximate stationarity of Problem~\eqref{prob:main} through its
KKT residual.
\begin{definition}
\label{def:eps-kkt}
Given \(\epsilon>0\), a triple \((\bx,\by,\bq)\) with
\(\bq\in\partial h(\by)\) is an \(\epsilon\)-KKT point of
Problem~\eqref{prob:main} if
\begin{equation}
    \mathcal{R}_{\mathrm{KKT}}(\bx,\by,\bq)
    :=
    \max\Bigl\{
        \bigl\|\nabla f(\bx)+J_{\calA}(\bx)^\top\bq\bigr\|,\;
        \|\calA(\bx)-\by\|
    \Bigr\} \leq \epsilon .
    \label{eq:kkt-residual}
\end{equation}
\end{definition}
At \(\epsilon=0\), \eqref{eq:kkt-residual} reduces to \(\by=\calA(\bx)\) and
\(\zero\in\nabla f(\bx)+J_{\calA}(\bx)^\top\partial h(\calA(\bx))\), the
first-order KKT condition of Problem~\eqref{prob:main}~\cite{li2020understanding}.

For complexity analysis, we specify the problem class~\cite{carmon2020lower,nemirovski1983problem}
and the oracle model for Problem~\eqref{prob:main}.
\begin{definition}
\label{def:problem-class}
Let \(\Delta>0\) and
\(\mathbf C=(L_f,G_f,L_{\calA},G_{\calA},L_h)\in\RR^5_{\ge0}\).
We denote by \(\mathcal I(\mathbf C,\Delta)\) the set of all triples
\((f,\calA,h)\), over all dimensions \(n,m\in\NN\), such that
Assumptions~\ref{ass:smoothness} and~\ref{ass:h} hold with the constants
in \(\mathbf C\) and
\(
    F(\zero)-\inf_{\bx\in\RR^n}F(\bx)\le\Delta .
\)
\end{definition}

\begin{definition}
\label{def:oracle-model}
We assume access to a deterministic first-order oracle
\[
    \mathcal{O}_{f,\mathcal{A}}(\bx)
    :=
    \bigl(f(\bx),\nabla f(\bx),\mathcal{A}(\bx),J_{\mathcal{A}}(\bx)\bigr),
\]
and an exact proximal oracle
\[
    \mathcal{O}_{h}(\bv,t)
    :=
    \operatorname{prox}_{t h}(\bv),
    \qquad t>0.
\]
Each call to either oracle counts as one oracle evaluation.
\end{definition}

\subsection{\texorpdfstring{Regularized augmented Lagrangian for Problem~\eqref{prob:main}}{Regularized augmented Lagrangian for Problem (P)}}
\label{subsec:regularized-al}
For \(\rho>0\), the augmented Lagrangian of Problem~\eqref{prob:main}
with multiplier \(\blmbd\in\RR^m\) is
\(\calL_\rho(\bx,\by,\blmbd):=f(\bx)+h(\by)+\langle\blmbd,\calA(\bx)-\by\rangle+\frac{\rho}{2}\|\calA(\bx)-\by\|^2\).
Eliminating \(\by\) and completing the square yields the reduced augmented
Lagrangian~\cite{rockafellar1976augmented,dhingra2019proximal},
\begin{equation}
    \calD_\rho(\bx,\blmbd)
    :=\min_{\by\in\RR^m}\calL_\rho(\bx,\by,\blmbd)
    \mathrel{=}f(\bx)+M_{1/\rho}^h\Bigl(\calA(\bx)+\frac{\blmbd}{\rho}\Bigr)
    -\frac{1}{2\rho}\|\blmbd\|^2,
    \label{eq:reduced-al}
\end{equation}
with unique minimizer and updated multiplier
\begin{equation}
    \by(\bx,\blmbd):=\prox_{h/\rho}\Bigl(\calA(\bx)+\frac{\blmbd}{\rho}\Bigr),
    \qquad
    \bq(\bx,\blmbd):=\blmbd+\rho\bigl[\calA(\bx)-\by(\bx,\blmbd)\bigr].
    \label{eq:prox-multiplier}
\end{equation}
The next lemma records the smoothness of \(\calD_\rho\) and shows that
maximizing it over \(\blmbd\) recovers \(F\).

\begin{lemma}
\label{lem:reduced-al}
Suppose that Assumptions~\ref{ass:smoothness} and~\ref{ass:h} hold, and let
\(\rho>0\).
\begin{enumerate}[label=(\alph*)]
\item \(\calD_\rho\) is continuously differentiable on \(\RR^n\times\RR^m\),
with \(\bq(\bx,\blmbd)\in\partial h(\by(\bx,\blmbd))\) and
\begin{equation}
    \nabla_{\bx}\calD_\rho(\bx,\blmbd)
    =\nabla f(\bx)+J_{\calA}(\bx)^\top\bq(\bx,\blmbd),
    \quad
    \nabla_{\blmbd}\calD_\rho(\bx,\blmbd)
    =\calA(\bx)-\by(\bx,\blmbd)
    =\frac{\bq(\bx,\blmbd)-\blmbd}{\rho}.
    \label{eq:reduced-al-gradient}
\end{equation}
For each \(\bx\in\RR^n\), \(\calD_\rho(\bx,\cdot)\) is concave and
\(\nabla_{\blmbd}\calD_\rho(\bx,\cdot)\) is \((1/\rho)\)-Lipschitz
continuous.
\item For every \((\bx,\blmbd)\in\RR^n\times\RR^m\),
\begin{equation}
    \calD_\rho(\bx,\blmbd)
    =\max_{\bq\in\RR^m}\Bigl\{\calL(\bx,\bq)-\frac{1}{2\rho}\|\bq-\blmbd\|^2\Bigr\},
    \label{eq:dual-smoothing}
\end{equation}
attained uniquely at \(\bq(\bx,\blmbd)\). Consequently,
\(\max_{\blmbd\in\RR^m}\calD_\rho(\bx,\blmbd)=F(\bx)\) for all \(\bx\in\RR^n\),
with maximizer set \(\partial h(\calA(\bx))\).
\end{enumerate}
\end{lemma}

By Lemma~\ref{lem:reduced-al}, Problem~\eqref{prob:main} reduces to the
unconstrained smooth NC-C minimax problem \(\min_{\bx\in\RR^n}\max_{\blmbd\in\RR^m}\calD_\rho(\bx,\blmbd)\).
By~\eqref{eq:kkt-residual} and~\eqref{eq:reduced-al-gradient},
\begin{equation}
    \mathcal R_{\mathrm{KKT}}\bigl(\bx,\by(\bx,\blmbd),\bq(\bx,\blmbd)\bigr)
    =\max\bigl\{
        \|\nabla_{\bx}\calD_\rho(\bx,\blmbd)\|,\;
        \|\nabla_{\blmbd}\calD_\rho(\bx,\blmbd)\|
    \bigr\},
    \label{eq:kkt-gs-identity}
\end{equation}
so any \(\epsilon\)-GS point
\((\bx,\blmbd)\)~\cite{pan2026optimal}, where the right-hand
side is at most \(\epsilon\), directly gives an \(\epsilon\)-KKT point
\((\bx,\by(\bx,\blmbd),\bq(\bx,\blmbd))\) of
Problem~\eqref{prob:main}; see also~\cite{li2025nonsmooth} for the
nonsmooth composite setting.

In this view, linearized ALM replaces exact \(\bx\)-minimization by one
gradient step and becomes alternating GDA on \(\calD_\rho\),
\begin{equation}
\left\{
\begin{aligned}
    \bx^{k+1}
    &=\bx^k-\eta\nabla_{\bx}\calD_\rho(\bx^k,\blmbd^k),\\
    \blmbd^{k+1}
    &=\blmbd^k+\rho\nabla_{\blmbd}\calD_\rho(\bx^{k+1},\blmbd^k)
    =\bq(\bx^{k+1},\blmbd^k),
\end{aligned}
\right.
\label{eq:alm-reduced}
\end{equation}
with \(\eta>0\). Since single-timescale GDA can fail to
converge even on bilinear problems~\cite{daskalakis2018training},
we adopt the primal~\cite{zhang2020smoothed} and
dual~\cite{lu2020hybrid,xu2023unified,zheng2023universal,li2026smoothing}
quadratic regularizers of smooth minimax methods and define the
regularized augmented Lagrangian for Problem~\eqref{prob:main},
\begin{equation}
    \mathcal P_s(\bx,\blmbd;\bz,\bmu)
    :=\calD_\rho(\bx,\blmbd)+\frac{1}{2\gamma}\|\bx-\bz\|^2
    -\frac{1}{2s}\|\blmbd-\bmu\|^2,
    \label{eq:potential}
\end{equation}
where \(\gamma,s>0\) and \((\bz,\bmu)\in\RR^n\times\RR^m\) are centers.
To stabilize the updates via strong convexity-strong concavity, our single-loop
algorithms alternate one primal-dual step on \(\mathcal P_s\) with center
updates for \((\bz,\bmu)\). Set \(L_{\calD,\bx}:=L_f+L_{\calA}L_h+\rho G_{\calA}^{2}\),
\(\mu_\gamma:=\gamma^{-1}-L_{\calD,\bx}\),
\(L_x:=\gamma^{-1}+L_{\calD,\bx}\), and \(L_\lambda:=\rho^{-1}+s^{-1}\).
The next proposition records the curvature, block-Lipschitz smoothness, and
saddle-point properties of \(\mathcal P_s\), which follow directly from
standard convex analysis~\cite{rockafellar1970convex}.
\begin{proposition}
\label{prop:curvature-smoothness}
Suppose that Assumptions~\ref{ass:smoothness} and~\ref{ass:h} hold, and let
\(\gamma>0\) satisfy \(\gamma L_{\calD,\bx}<1/2\). Then, for every \(s>0\) and
\((\bz,\bmu)\in\RR^n\times\RR^m\), the following statements hold.
\begin{enumerate}[label=(\alph*)]
\item \(\mathcal P_s(\cdot,\cdot;\bz,\bmu)\) is \(\mu_\gamma\)-strongly
convex in \(\bx\) and \(s^{-1}\)-strongly concave in \(\blmbd\). Moreover, for
all \(\bx,\bx'\in\RR^n\) and \(\blmbd,\blmbd'\in\RR^m\),
\begin{align*}
    \bigl\|
    \nabla_{\bx}\mathcal{P}_s(\bx,\blmbd;\bz,\bmu)
    -\nabla_{\bx}\mathcal{P}_s(\bx',\blmbd';\bz,\bmu)
    \bigr\|
    &\le L_x\|\bx-\bx'\|+G_{\calA}\|\blmbd-\blmbd'\|,\\
    \bigl\|
    \nabla_{\blmbd}\mathcal{P}_s(\bx,\blmbd;\bz,\bmu)
    -\nabla_{\blmbd}\mathcal{P}_s(\bx',\blmbd';\bz,\bmu)
    \bigr\|
    &\le G_{\calA}\|\bx-\bx'\|+L_\lambda\|\blmbd-\blmbd'\|.
\end{align*}
In particular, \(\mathcal P_s(\cdot,\cdot;\bz,\bmu)\) is \(L_x\)-smooth in
\(\bx\) and \(L_\lambda\)-smooth in \(\blmbd\).
\item \(\mathcal P_s(\cdot,\cdot;\bz,\bmu)\) admits a unique saddle point
\((\bx^\star,\blmbd^\star)\in\RR^n\times\RR^m\) and satisfies
\[
    \min_{\bx\in\RR^n}\max_{\blmbd\in\RR^m}\mathcal P_s(\bx,\blmbd;\bz,\bmu)
    =\max_{\blmbd\in\RR^m}\min_{\bx\in\RR^n}\mathcal P_s(\bx,\blmbd;\bz,\bmu)
    =\mathcal P_s(\bx^\star,\blmbd^\star;\bz,\bmu).
\]
\end{enumerate}
\end{proposition}

\section{Tight Complexity of the LRALM Family}
\label{sec:algorithm-framework}

In this section, we address \hyperlink{Q1}{\textbf{(Q1)}} by establishing tight
iteration-complexity bounds for single-loop linearized ALMs on
Problem~\eqref{prob:main}. Taking alternating GDA steps on \(\mathcal P_s\)
in~\eqref{eq:potential} together with moving primal and dual centers gives the
LRALM method in Algorithm~\ref{alg:lralm}.

\begin{definition}
\label{def:splalm}
The method \(\mathrm{LRALM}_{\theta,\alpha,\beta}\) is defined by
Algorithm~\ref{alg:lralm} with output
\((\bx^{k+1},\by^{k+1},\bq^k)
:=(\bx^{k+1},\by(\bx^{k+1},\blmbd^k),\bq(\bx^{k+1},\blmbd^k))\), where
\(\by(\cdot,\cdot)\) and \(\bq(\cdot,\cdot)\) are defined
in~\eqref{eq:prox-multiplier}.
\end{definition}

The LRALM family \(\calC_\theta\), formally defined in
Definition~\ref{def:complexity-measure}, collects all such methods over
\((\alpha,\beta)\in[0,1]^2\).

\begin{algorithm}[t]
    \caption{\(\mathrm{LRALM}_{\theta,\alpha,\beta}\), a member of \(\calC_\theta\)}
    \label{alg:lralm}
    \small
    \begin{algorithmic}[1]
        \STATE \textbf{Input:} Parameters \(\theta=(\rho,s,\gamma,\eta,\tau)\) with
        \(\rho>0\), \(0<\eta\le\gamma\), and \(0<\tau\le s\);
        center stepsizes \((\alpha,\beta)\in[0,1]^2\); iteration limit \(K\).
        \STATE \textbf{Initialize:} \((\bx^0,\blmbd^0,\bz^0,\bmu^0)\in\RR^n\times\RR^m\times\RR^n\times\RR^m\).
        \FOR{\(k=0,1,\ldots,K-1\)}
            \STATE \(\bx^{k+1}=\bx^k-\eta\Bigl(\nabla_{\bx}\calL_\rho\bigl(\bx^k,\by(\bx^k,\blmbd^k),\blmbd^k\bigr)+\tfrac{1}{\gamma}(\bx^k-\bz^k)\Bigr)\).
            \STATE \(\blmbd^{k+1}=\blmbd^k+\tau\Bigl(\calA(\bx^{k+1})-\by(\bx^{k+1},\blmbd^k)-\tfrac{1}{s}(\blmbd^k-\bmu^k)\Bigr)\).
            \STATE \(\bz^{k+1}=\bz^k+\alpha(\bx^{k+1}-\bz^k)\).
            \STATE \(\bmu^{k+1}=\bmu^k+\beta(\blmbd^{k+1}-\bmu^k)\).
        \ENDFOR
        \STATE \textbf{Output:}
        \(\bigl\{\bigl(\bx^{k+1},\,\by(\bx^{k+1},\blmbd^k),\,\bq(\bx^{k+1},\blmbd^k)\bigr)\bigr\}_{k=0}^{K-1}\).
    \end{algorithmic}
\end{algorithm}

When applied to \(\calD_\rho\), existing single-loop methods from both ALM
and smooth minimax optimization take the form of Algorithm~\ref{alg:lralm},
with different choices of the dual regularizer \(s\) and the center
stepsizes \((\alpha,\beta)\). We first clarify the effect of \(s\) on
\(\calD_\rho\) and then organize these center-movement choices into the LRALM
family.

\subsection{The role of the dual regularizer and the LRALM family}
\label{subsec:lralm-family}

In linearized ALM~\eqref{eq:alm-reduced} on \(\calD_\rho\),
constraint violation is typically controlled by choosing a large penalty
\(\rho\)~\cite{luSingleLoopGradientDescent2022,sun2024dual_descent_alm_admm,elbourkhissi2025perturbed}.
Indeed, combining~\eqref{eq:reduced-al-gradient} with the bound
\(\|\bq(\bx,\blmbd)\|\le L_h\) from Lemma~\ref{lem:reduced-al}(a) and
Assumption~\ref{ass:h} gives
\[
    \|\calA(\bx)-\by(\bx,\blmbd)\|
    \le\frac{L_h+\|\blmbd\|}{\rho}.
\]
On the other hand, since \(L_{\calD,\bx}\) grows linearly with \(\rho\),
a large \(\rho\) increases the primal smoothness modulus in
Proposition~\ref{prop:curvature-smoothness}(a) and restricts the primal stepsize
\(\eta\). The next lemma shows that,
without increasing \(\rho\), maximizing over \(\blmbd\) with the dual regularizer
raises the effective penalty from \(\rho\) to \(\rho+s\). The identity
follows from~\eqref{eq:dual-smoothing}, since the infimal convolution of
\(\|\cdot\|^2/(2\rho)\) and \(\|\cdot\|^2/(2s)\) is
\(\|\cdot\|^2/(2(\rho+s))\); the full proof is given in
the supplementary material.

\begin{lemma}
\label{lem:effective-penalty}
Suppose that Assumption~\ref{ass:h} holds, and let \(\rho,s>0\). For every
\((\bx,\bmu)\in\RR^n\times\RR^m\),
\begin{equation}
    \max_{\blmbd\in\RR^m}\Bigl\{\calD_\rho(\bx,\blmbd)
    -\frac{1}{2s}\|\blmbd-\bmu\|^2\Bigr\}
    =\calD_{\rho+s}(\bx,\bmu).
    \label{eq:effective-penalty}
\end{equation}
Moreover, the left-hand side of~\eqref{eq:effective-penalty} has a unique
maximizer \(\blmbd^\star(\bx,\bmu)\), which satisfies
\[
    \by\bigl(\bx,\blmbd^\star(\bx,\bmu)\bigr)
    =\prox_{h/(\rho+s)}\bigl(\calA(\bx)+\bmu/(\rho+s)\bigr),
    \qquad
    \bq^\star:=\bq\bigl(\bx,\blmbd^\star(\bx,\bmu)\bigr)
    \in\partial h\bigl(\by(\bx,\blmbd^\star(\bx,\bmu))\bigr),
\]
and
\begin{equation}
\begin{aligned}
    \blmbd^\star(\bx,\bmu)
    &=\frac{s}{\rho+s}\bq^\star+\frac{\rho}{\rho+s}\bmu,\\
    \nabla_{\bmu}\calD_{\rho+s}(\bx,\bmu)
    &=\calA(\bx)-\by\bigl(\bx,\blmbd^\star(\bx,\bmu)\bigr)
    =\frac{\bq^\star-\bmu}{\rho+s}
    =\frac{\blmbd^\star(\bx,\bmu)-\bmu}{s}.
\end{aligned}
    \label{eq:effective-penalty-maximizer}
\end{equation}
\end{lemma}

Under Assumption~\ref{ass:h}, we have \(\|\bq^\star\|\le L_h\),
so~\eqref{eq:effective-penalty-maximizer} replaces \(\rho\) in the previous
bound by \(\rho+s\) and gives
\[
    \bigl\|\calA(\bx)-\by\bigl(\bx,\blmbd^\star(\bx,\bmu)\bigr)\bigr\|
    \le\frac{L_h+\|\bmu\|}{\rho+s}.
\]
Meanwhile, the primal smoothness modulus \(L_x=\gamma^{-1}+L_{\calD,\bx}\)
of \(\mathcal P_s\) is independent of \(s\) by
Proposition~\ref{prop:curvature-smoothness}(a). Increasing \(s\) therefore
strengthens the effective penalty without shrinking the primal stepsize
\(\eta\).

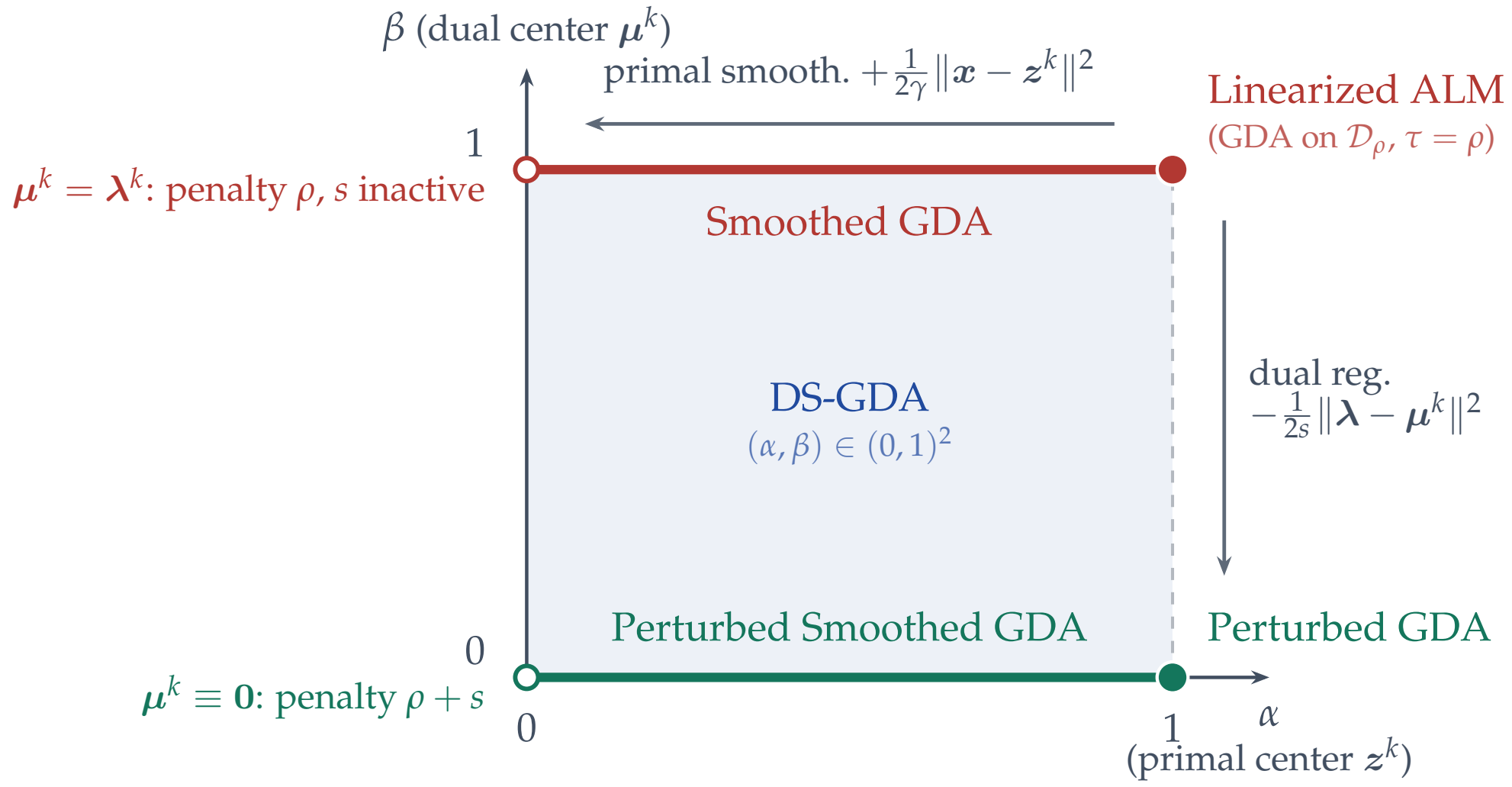
\begin{figure}[t]
\centering
\definecolor{cSlate}{RGB}{65,78,96}
\definecolor{cRed}{RGB}{178,56,50}
\definecolor{cTeal}{RGB}{20,118,92}
\definecolor{cNavy}{RGB}{32,74,158}
\resizebox{\ifdim\width>\textwidth\textwidth\else\width\fi}{!}{%
\begin{tikzpicture}[x=5.6cm,y=4.4cm,font=\small,
  axis/.style={cSlate,line width=0.85pt,-{Stealth[length=4.5pt]}},
  mod/.style={cSlate!85,line width=0.9pt,-{Stealth[length=4.8pt]}}]
  \fill[cNavy!8] (0,0) rectangle (1,1);
  \draw[cSlate!40,dashed,line width=0.75pt] (1,0)--(1,1);
  \draw[axis] (0,-0.02)--(0,1.20) node[above,text=cSlate] {$\beta$\ \footnotesize(dual center $\bmu^k$)};
  \draw[axis] (-0.02,0)--(1.15,0) node[below=3pt,text=cSlate,align=center] {$\alpha$\\[-1pt]\footnotesize(primal center $\bz^k$)};
  \draw[cRed,line width=2.3pt] (0,1)--(1,1);
  \draw[cTeal,line width=2.3pt] (0,0)--(1,0);
  \node[circle,draw=cRed,fill=white,line width=1.1pt,inner sep=2.0pt] at (0,1) {};
  \node[circle,draw=cTeal,fill=white,line width=1.1pt,inner sep=2.0pt] at (0,0) {};
  \node[circle,draw=white,fill=cRed,line width=0.8pt,inner sep=2.7pt] at (1,1) {};
  \node[circle,draw=white,fill=cTeal,line width=0.8pt,inner sep=2.7pt] at (1,0) {};
  \node[below,text=cSlate] at (0,-0.04) {$0$};
  \node[below,text=cSlate] at (1,-0.04) {$1$};
  \node[left,align=right] at (-0.04,1) {\textcolor{cSlate}{$1$}\\[-1pt]
    \footnotesize\textcolor{cRed}{$\bmu^k=\blmbd^k$: penalty $\rho$, $s$ inactive}};
  \node[left,align=right] at (-0.04,0) {\textcolor{cSlate}{$0$}\\[-1pt]
    \footnotesize\textcolor{cTeal}{$\bmu^k\equiv\zero$: penalty $\rho+s$}};
  \node[anchor=north,text=cRed] at (0.5,0.96)
    {Smoothed GDA};
  \node[align=center,text=cNavy] at (0.5,0.5)
    {DS-GDA\\[-1pt]
    \scriptsize\textcolor{cNavy!75}{$(\alpha,\beta)\in(0,1)^2$}};
  \node[anchor=south,text=cTeal] at (0.5,0.04)
    {Perturbed Smoothed GDA};
  \node[cRed,anchor=south west,align=left] at (1.03,0.99)
    {Linearized ALM\\[-1pt]
    \scriptsize\textcolor{cRed!80}{(GDA on $\calD_\rho$, $\tau=\rho$)}};
  \node[cTeal,anchor=south west] at (1.03,0.04)
    {Perturbed GDA};
  \draw[mod] (0.91,1.09)--(0.09,1.09)
     node[midway,above=1.5pt,font=\footnotesize,text=cSlate]
     {primal smooth.\ $+\frac{1}{2\gamma}\|\bx-\bz^k\|^2$};
  \draw[mod] (1.08,0.90)--(1.08,0.20)
     node[midway,right=2pt,align=left,font=\footnotesize,text=cSlate]
     {dual reg.\\[-1pt]$-\frac{1}{2s}\|\blmbd-\bmu^k\|^2$};
\end{tikzpicture}%
}
\caption{Center-stepsize space \((\alpha,\beta)\in[0,1]^2\) of the LRALM family
\(\calC_\theta\) acting on \(\calD_\rho\). Red, teal, and blue mark the top
edge \(\beta=1\), the bottom edge \(\beta=0\) (with \(\bmu^0=\zero\)), and the
region \((0,1)^2\), while solid circles mark the corner methods at
\(\alpha=1\) and open circles mark the excluded endpoints at \(\alpha=0\).}
\label{fig:lralm-family}
\end{figure}

In Algorithm~\ref{alg:lralm}, however, the dual center \(\bmu^k\) tracks
\(\blmbd^k\) with stepsize \(\beta\), which determines whether this effective
penalty persists across iterations. At \(\beta=1\), we have
\(\bmu^k=\blmbd^k\) for all \(k\ge1\), so the dual regularizer vanishes and the
effective penalty reverts to \(\rho\), whereas at \(\beta=0\) with
\(\bmu^0=\zero\) the dual center stays fixed at \(\zero\) and retains the full
\(\rho+s\) of Lemma~\ref{lem:effective-penalty}. Varying the constant center
stepsizes \((\alpha,\beta)\in[0,1]^2\) yields the LRALM family \(\calC_\theta\),
formalized in Definition~\ref{def:complexity-measure} and illustrated in
Figure~\ref{fig:lralm-family}.

\begin{definition}
\label{def:complexity-measure}
Given \(\theta=(\rho,s,\gamma,\eta,\tau)\) as in Definition~\ref{def:splalm},
the LRALM family is
\(\calC_\theta:=\{\mathrm{LRALM}_{\theta,\alpha,\beta}:(\alpha,\beta)\in[0,1]^2\}\).
Given \(\mathbf C\in\RR^5_{\ge0}\), define the admissible initialization set
\[
    \mathcal S_0(\mathbf C)
    :=\bigl\{\bu^0=(\bx^0,\blmbd^0,\bz^0,\bmu^0):\
        \bx^0=\zero,\ \bmu^0=\zero,\ \|\blmbd^0\|\le L_h,\ \|\bz^0\|\le1\bigr\}.
\]
For \(\epsilon>0\), \(I\in\mathcal I(\mathbf C,\Delta)\), and
\(\bu^0\in\mathcal S_0(\mathbf C)\), the iteration complexity of
\(\mathrm{LRALM}_{\theta,\alpha,\beta}\) on \((I,\bu^0)\) is
\[
    T_\epsilon(\theta,\alpha,\beta;I,\bu^0)
    :=\inf\Bigl\{K\ge1:\ \min_{0\le k<K}
        \mathcal R_{\mathrm{KKT}}(\bx^{k+1},\by^{k+1},\bq^k)\le\epsilon\Bigr\},
\]
with \(\inf\varnothing=+\infty\). The worst-case iteration complexity of
\(\mathrm{LRALM}_{\theta,\alpha,\beta}\) over \(\mathcal I(\mathbf C,\Delta)\) is
\[
    T_\epsilon\bigl(\theta,\alpha,\beta;\mathcal I(\mathbf C,\Delta)\bigr)
    :=\sup_{\substack{I\in\mathcal I(\mathbf C,\Delta)\\ \bu^0\in\mathcal S_0(\mathbf C)}}\
    T_\epsilon(\theta,\alpha,\beta;I,\bu^0),
\]
and that of \(\calC_\theta\) is
\[
    T_\epsilon\bigl(\calC_\theta,\mathcal I(\mathbf C,\Delta)\bigr)
    :=\inf_{(\alpha,\beta)\in[0,1]^2}\
    T_\epsilon\bigl(\theta,\alpha,\beta;\mathcal I(\mathbf C,\Delta)\bigr).
\]
\end{definition}

Figure~\ref{fig:lralm-family} locates existing single-loop methods applied
to \(\calD_\rho\) in the \((\alpha,\beta)\) plane of \(\calC_\theta\). On the
top edge \(\beta=1\), \(\alpha=1\) with \(\tau=\rho\) recovers linearized
ALM~\eqref{eq:alm-reduced} and \(\alpha\in(0,1)\) gives
Smoothed GDA~\cite{zhang2020smoothed}; on the bottom edge \(\beta=0\)
with \(\bmu^0=\zero\), \(\alpha=1\) yields Perturbed
GDA~\cite{lu2020hybrid,xu2023unified} and \(\alpha\in(0,1)\) yields
Perturbed Smoothed GDA~\cite{li2026smoothing}. The interior
\((\alpha,\beta)\in(0,1)^2\) corresponds to doubly smoothed GDA
(DS-GDA)~\cite{zheng2023universal}.

\subsection{\texorpdfstring{Tight bounds in two effective-penalty regimes}{Tight bounds in two effective-penalty regimes}}
\label{sec:lower-bound}

We now establish matching lower and upper bounds for \(\calC_\theta\) in
two regimes of the effective penalty \(\rho+s\); the proofs are given in

the supplementary material. By reusing
\(\mathcal O_{f,\calA}(\bx^{k+1})\) in the next primal step, each iteration
of Algorithm~\ref{alg:lralm} makes one call to \(\mathcal O_{f,\calA}\) and
two calls to \(\mathcal O_h\), so the iteration bounds below agree with the
oracle complexity up to a constant factor.

\subsubsection{Tight bounds for \texorpdfstring{\(\rho+s\le c_1/\epsilon\)}{rho+s <= c1/epsilon}}

We first establish that, while \(\epsilon(\rho+s)\) stays bounded by a
constant, \(\calC_\theta\) attains \(\epsilon^{-4}\) and no
center-stepsize pair \((\alpha,\beta)\in[0,1]^2\) improves on it.

\begin{theorem}
\label{thm:lower-bound}
Let \(\mathbf C\in\RR^5_{\ge0}\) satisfy \(L_f,G_f,G_{\calA},L_h\ge1\),
and let \(0<\eta\le\gamma\), \(\tau>0\), and \(\Delta>0\). Then there exist
constants \(c_1,c_2>0\), independent of \((\rho,s)\), such
that, for every \(\epsilon\in(0,c_2)\), every \(\rho>0\), and every
\(s\ge\tau\) with \(\epsilon(\rho+s)\le c_1\),
\[
    T_\epsilon\bigl(\calC_\theta,\mathcal I(\mathbf C,\Delta)\bigr)
    =\Omega(\epsilon^{-4}),
\]
where the constant in \(\Omega(\cdot)\) is independent of \((\rho,s)\).
\end{theorem}
The matching upper bound is attained by Smoothed GDA~\cite{zhang2020smoothed},
the member of \(\calC_\theta\) with \(\beta=1\). Since
\(\bmu^k=\blmbd^k\) for all \(k\ge1\) at \(\beta=1\), the LRALM updates
from \(k=1\) on coincide with Smoothed GDA and no longer involve \(s\);
the multipliers also remain in the ball of radius \(L_h\), so the projection
onto the compact dual domain is inactive.
Consequently, the bound follows
from~\cite{li2026smoothing} with \(\alpha=c\,\epsilon^2\); the
full proof is given in the supplementary material.

\begin{proposition}
\label{prop:upper-beta1}
Let \(\mathbf C\in\RR^5_{\ge0}\), \(\Delta>0\), and \(\rho>0\), set
\(\ell:=\max\{L_{\calD,\bx},G_{\calA},1/\rho\}\), and let
\(0<\gamma\le\bar c/\ell\), \(0<\eta\le\gamma/2\), and
\(0<\tau\le\bar c(\eta/\gamma)^2/\ell\) for a sufficiently small
constant \(\bar c>0\). Then there exists \(c>0\) such that, for
every \(s\ge\tau\) and every \(\epsilon\in(0,1]\), the choice
\((\alpha,\beta)=(c\,\epsilon^2,1)\) yields
\[
    T_\epsilon\bigl(\theta,c\,\epsilon^2,1;\mathcal I(\mathbf C,\Delta)\bigr)
    =\mathcal O(\epsilon^{-4}),
\]
where the constant in \(\mathcal O(\cdot)\) is independent of \(s\).
\end{proposition}
\noindent Under the stepsize conditions of
Proposition~\ref{prop:upper-beta1}, combining Theorem~\ref{thm:lower-bound}
and Proposition~\ref{prop:upper-beta1}
yields \(T_\epsilon(\calC_\theta,\mathcal I(\mathbf C,\Delta))=\Theta(\epsilon^{-4})\)
whenever \(\rho+s\le c_1/\epsilon\). Since \(c_1\) and \(c_2\) are
independent of \((\rho,s)\), improving on \(\Theta(\epsilon^{-4})\) over
\(\mathcal I(\mathbf C,\Delta)\) requires \(\rho+s>c_1/\epsilon\) for
all sufficiently small \(\epsilon\), which answers the first part
of~\hyperlink{Q1}{\textbf{(Q1)}}. Specifically, with a constant penalty
\(\rho\) and a dual regularization parameter \(s\) that does not grow as
\(\epsilon\to0\), we have \(\rho+s\le c_1/\epsilon\) for all sufficiently small
\(\epsilon\), so \(\calC_\theta\) cannot improve on \(\Theta(\epsilon^{-4})\)
and hence does not attain \(\mathcal O(\epsilon^{-3})\).

\subsubsection{Tight bounds for \texorpdfstring{\(s=c_s/\epsilon\)}{s = cs/epsilon}}

To reach the regime \(\rho+s>c_1/\epsilon\) with a constant \(\rho>0\), we
scale \(s=c_s/\epsilon\), transferring the \(\epsilon^{-1}\) growth to the
dual regularizer while preserving the primal stepsize. At this scale,
\(\calC_{\theta(\epsilon)}\)
attains \(\epsilon^{-3}\), and no center-stepsize pair
\((\alpha,\beta)\in[0,1]^2\) improves on it.

\begin{theorem}
\label{thm:lower-bound-fixed-rho}
Let \(\mathbf C\in\RR^5_{\ge0}\) satisfy \(L_f,G_f,G_{\calA},L_h\ge1\),
and let \(\rho>0\), \(0<\eta\le\gamma\), \(\tau>0\), \(c_s>0\), and
\(\Delta>0\). Then there exist constants \(B\ge1\) and
\(\epsilon_0\in(0,1]\) such that, for every
\(\epsilon<\epsilon_0\) and every \((\alpha,\beta)\in[0,1]^2\), with
\(\theta(\epsilon):=(\rho,c_s/\epsilon,\gamma,\eta,\tau)\),
\[
    T_\epsilon\bigl(\theta(\epsilon),\alpha,\beta;\mathcal I(\mathbf C,\Delta)\bigr)
    \begin{cases}
        =+\infty, & \alpha=0\ \text{or}\ \alpha>B\epsilon,\\[2pt]
        =\Omega\bigl(\alpha^{-1}\epsilon^{-2}\bigr), & 0<\alpha\le B\epsilon.
    \end{cases}
\]
In particular,
\(T_\epsilon(\calC_{\theta(\epsilon)},\mathcal I(\mathbf C,\Delta))=\Omega(\epsilon^{-3})\).
\end{theorem}
\noindent The matching upper bound is attained by Perturbed Smoothed
GDA, the member of \(\calC_{\theta(\epsilon)}\) with
\(\beta=0\). At \(\beta=0\), \(\bmu^k\equiv\zero\) and LRALM coincides
with Perturbed Smoothed GDA with dual perturbation \(r_y=1/s\). A sufficiently
large \(c_s\) keeps the induced stationarity bias, of order \(1/s\), at most
\(\epsilon/4\), and the analysis of~\cite{li2026smoothing} admits
\(\alpha=c_\alpha\epsilon\) for a sufficiently small \(c_\alpha>0\) and yields
the bound; the full proof
is given in the supplementary material.

\begin{proposition}
\label{prop:upper-beta0}
Let \(\mathbf C\in\RR^5_{\ge0}\), \(\Delta>0\), and \(\rho>0\), let
\(\ell,\gamma,\eta,\tau\) be as in Proposition~\ref{prop:upper-beta1}, and let
\(c_s>0\) be sufficiently large. Then there exist constants \(c_\alpha>0\)
and \(\bar\epsilon\in(0,1]\) such that, for every \(\epsilon\le\bar\epsilon\), with
\(\theta(\epsilon):=(\rho,c_s/\epsilon,\gamma,\eta,\tau)\), the choice
\((\alpha,\beta)=(c_\alpha\epsilon,0)\) yields
\[
    T_\epsilon\bigl(\theta(\epsilon),c_\alpha\epsilon,0;\mathcal I(\mathbf C,\Delta)\bigr)
    =\mathcal O(\epsilon^{-3}).
\]
\end{proposition}
\noindent Combining Theorem~\ref{thm:lower-bound-fixed-rho} and
Proposition~\ref{prop:upper-beta0} yields the tight family complexity
\(T_\epsilon(\calC_{\theta(\epsilon)},\mathcal I(\mathbf C,\Delta))
=\Theta(\epsilon^{-3})\) for every sufficiently large \(c_s\).

Taken together, the two regimes answer~\hyperlink{Q1}{\textbf{(Q1)}}.
With a constant penalty, growing \(s\) as \(c_s/\epsilon\) with a
sufficiently large \(c_s\) reduces the family complexity from
\(\Theta(\epsilon^{-4})\) to \(\Theta(\epsilon^{-3})\), while keeping the dual
stepsize \(\tau\) independent of \(\epsilon\) as well. Some
\(\epsilon\)-dependent growth in \(s\) is necessary, since accuracy-independent
\((\rho,s)\) yield only \(\Theta(\epsilon^{-4})\) by
Theorem~\ref{thm:lower-bound}. Within
\(\calC_{\theta(\epsilon)}\), however, growing \(s\) alone does not improve
on \(\epsilon^{-3}\), since Theorem~\ref{thm:lower-bound-fixed-rho}
holds for every \(c_s>0\) and every \((\alpha,\beta)\in[0,1]^2\). The
limitation comes from the primal center stepsize \(\alpha\), that is, from how
the primal center moves, independently of the dual center. This motivates
the accelerated method in Section~\ref{sec:acceleration}.

\begin{remark}\label{rem:lralm-positioning}
\begin{enumerate}[label=(\roman*)]
\item Recent lower bounds establish the optimal first-order complexity for
general smooth NC-C minimax problems~\cite{pan2026optimal,zhang2026matching}
and certify the tightness of individual single-loop methods under prescribed
stepsize
conditions~\cite{mahdavinia2022tight,li2026smoothing}.
Complementing these results, Theorems~\ref{thm:lower-bound}
and~\ref{thm:lower-bound-fixed-rho} establish lower bounds on the structured
minimax problems \(\min_{\bx}\max_{\blmbd}\calD_\rho(\bx,\blmbd)\)
generated by \(\mathcal I(\mathbf C,\Delta)\) across every constant
center-stepsize pair \((\alpha,\beta)\in[0,1]^2\). Specifically, because
each such instance is a smooth NC-C minimax problem, these bounds transfer
directly to the corresponding single-loop methods on
general smooth NC-C problems. In particular, Smoothed GDA and
Perturbed Smoothed GDA attain the tight worst-case orders within
\(\calC_\theta\) in their respective regimes.
\item In nonconvex constrained optimization, existing single-loop ALM analyses
typically establish complexity guarantees through growing
penalties~\cite{luSingleLoopGradientDescent2022,sun2024dual_descent_alm_admm,elbourkhissi2025perturbed},
while constant-penalty methods rely on uniform linear independence
constraint qualification~\cite{liu2026fixed} or exact prox-linear
subproblems~\cite{zhu2026nonconvex}.
Alongside these developments, combining Theorems~\ref{thm:lower-bound}
and~\ref{thm:lower-bound-fixed-rho} with Propositions~\ref{prop:upper-beta1}
and~\ref{prop:upper-beta0} establishes tight worst-case complexity bounds for
the constant-penalty LRALM family \(\calC_\theta\) on
Problem~\eqref{prob:main}.
While the open question of~\cite{Alacaoglu24a} targets deterministic
nonlinearly constrained problems, these bounds show that the family attains
\(\Theta(\epsilon^{-3})\) in this composite setting with a constant penalty and
\(\epsilon\)-dependent dual regularization.
\end{enumerate}
\end{remark}

\subsection{Proof sketch of the lower bounds}
\label{sec:lb-proof}
We sketch the proofs of Theorems~\ref{thm:lower-bound}
and~\ref{thm:lower-bound-fixed-rho} in two steps; the explicit constructions of
Instances~II and~IV and all proofs are given in the supplementary material.
Throughout this subsection, we write
\(R_k:=\mathcal R_{\mathrm{KKT}}(\bx^{k+1},\by^{k+1},\bq^k)\) and assume
\(L_f,G_f,G_{\calA},L_h\ge1\).

\textit{Step~1.} Following~\cite{mahdavinia2022tight}, we
partition \([0,1]^2\) by the primal center stepsize
\(\alpha\), with threshold \(\alpha_\epsilon=\Theta(\epsilon^2)\) for
\(\rho+s\le c_1/\epsilon\) and \(B\epsilon\) for \(s=c_s/\epsilon\). Below
either threshold, the plateau Instance~I forces
\(\Omega(\Delta/(\gamma\alpha\epsilon^2))\) iterations for every \(\beta\).

\textit{Step~2.} Above these thresholds, Instances~II, III, and~IV yield
\(T_\epsilon=+\infty\). Since \(\mathcal I(\mathbf C,\Delta)\) requires
\(f\) to be globally Lipschitz, these instances keep the iterates bounded and
the KKT residual from decaying, whereas the large-stepsize regime
of~\cite{mahdavinia2022tight} relies on divergence. Under
\(\rho+s\le c_1/\epsilon\), the rotating
Instance~II covers \(\beta\in(0,1]\) and the Huber Instance~III sustains a dual
feasibility bias at \(\beta=0\); under \(s=c_s/\epsilon\), the radial
Instance~IV covers every \(\beta\in[0,1]\). Taking the infimum over
\((\alpha,\beta)\in[0,1]^2\) of the bounds in the two steps yields both
theorems.

\subsubsection{Hard instances for \texorpdfstring{\(\rho+s\le c_1/\epsilon\)}{rho+s <= c1/epsilon}}

\paragraph{I. Plateau instance.}
The first instance covers \(\alpha\le\alpha_\epsilon\) for every
\(\beta\in[0,1]\), and is reused below for \(\alpha\le B\epsilon\).

\begin{example}
\label{ex:plateau}
Let \(0<\epsilon\le\min\{1/2,\sqrt{\Delta}/2\}\) and
\(x_\star:=\Delta/(2\epsilon)+\epsilon\). Define \(\calA\equiv0\),
\(h(y)=|y|\), and \(f\colon\RR\to\RR\) by
\[
f(x)=
\begin{cases}
2\epsilon(x_\star-x)-2\epsilon^2, & x<x_\star-2\epsilon,\\[4pt]
\dfrac12(x-x_\star)^2, & |x-x_\star|\le2\epsilon,\\[4pt]
2\epsilon(x-x_\star)-2\epsilon^2, & x>x_\star+2\epsilon.
\end{cases}
\]
\end{example}

Because \(\calA\equiv0\), the dual variables remain zero from \(\bu^0=\zero\),
and every iterate on \(x<x_\star-2\epsilon\) satisfies
\(R_k=|\nabla f(x^{k+1})|=2\epsilon\) and \(|x^k-z^k|\le2\gamma\epsilon\). A
small \(\alpha\) therefore limits each center step to order
\(\gamma\alpha\epsilon\), yielding the lower bound below; its proof is given in
the supplementary material.

\begin{lemma}
\label{lem:main-plateau}
Let \((\alpha,\beta)\in[0,1]^2\), and let \(\epsilon>0\) satisfy the
condition of Example~\ref{ex:plateau} and \(\epsilon^2\le\Delta/(16\gamma)\).
Then the instance \(I\) of Example~\ref{ex:plateau} belongs to
\(\mathcal I(\mathbf C,\Delta)\), and \(\mathrm{LRALM}_{\theta,\alpha,\beta}\)
initialized at \(\bu^0=\zero\) satisfies
\(T_\epsilon(\theta,\alpha,\beta;I,\zero)>\Delta/(16\gamma\alpha\epsilon^2)\)
for \(\alpha>0\), and \(T_\epsilon(\theta,\alpha,\beta;I,\zero)=+\infty\) for
\(\alpha=0\).
\end{lemma}

\paragraph{II. Rotating instance.}
For \(\alpha>\alpha_\epsilon\) and \(\beta\in(0,1]\), we build on the
quadratic coupling of~\cite{li2026smoothing} with a radial
extension that induces a non-decaying unit-modulus mode in the LRALM
recursion, covering every \(\alpha>\alpha_\epsilon\); the
explicit instance and proof are given in the supplementary material.

\begin{lemma}
\label{lem:main-rotation}
Under the assumptions of Theorem~\ref{thm:lower-bound}, let
\(\alpha>\alpha_\epsilon\) and \(\beta\in(0,1]\). Then there exists an
instance \(I\in\mathcal I(\mathbf C,\Delta)\) on which
\(\mathrm{LRALM}_{\theta,\alpha,\beta}\) initialized at \(\bu^0=\zero\)
satisfies \(R_k\ge\frac76\epsilon\) for all \(k\ge0\); hence
\(T_\epsilon(\theta,\alpha,\beta;I,\zero)=+\infty\).
\end{lemma}

\paragraph{III. Huber instance.}
For \(\alpha>\alpha_\epsilon\) and \(\beta=0\), we retain \(f\) from
Example~\ref{ex:plateau} and pair a constant mapping \(\calA\) with a Huber
penalty \(h\).

\begin{example}
\label{ex:huber}
Let \(\epsilon\) and \(f\) be as in Example~\ref{ex:plateau}, and set
\[
q_h:=2s+1,
\qquad
a_0:=\frac{4\epsilon(\rho+3s+1)}{q_h}.
\]
Define \(\calA\equiv a_0\) and \(h\colon\RR\to\RR\) by
\[
h(y)=
\begin{cases}
\dfrac{q_h}{2}y^2, & |y|\le q_h^{-1},\\[4pt]
|y|-\dfrac{1}{2q_h}, & |y|>q_h^{-1}.
\end{cases}
\]
\end{example}

At \(\beta=0\), the dual center remains at zero, so the regularizer pulls
\(\lambda^k\) below the multiplier needed for feasibility, sustaining a
dual bias. The proof is deferred to the supplementary material.

\begin{lemma}
\label{lem:main-huber}
Let \(\alpha\in[0,1]\) and \(\beta=0\), and suppose \(\epsilon\) satisfies
the condition of Example~\ref{ex:plateau} and \(\epsilon(\rho+s)\le c_1\).
Then the instance \(I\) of Example~\ref{ex:huber} belongs to
\(\mathcal I(\mathbf C,\Delta)\), and \(\mathrm{LRALM}_{\theta,\alpha,0}\)
initialized at \(\bu^0=\zero\) satisfies \(R_k\ge3\epsilon\) for all
\(k\ge0\); hence \(T_\epsilon(\theta,\alpha,0;I,\zero)=+\infty\).
\end{lemma}

Combining Lemmas~\ref{lem:main-plateau}--\ref{lem:main-huber} over the three
regions of \([0,1]^2\) yields Theorem~\ref{thm:lower-bound}, whose complete
proof appears in the supplementary material.

\subsubsection{Hard instance for \texorpdfstring{\(s=c_s/\epsilon\)}{s = cs/epsilon}}

For \(\alpha\le B\epsilon\) and any \(\beta\in[0,1]\),
Lemma~\ref{lem:main-plateau} on Example~\ref{ex:plateau} directly gives the
lower bound \(\Omega(\Delta\gamma^{-1}\epsilon^{-3})\).

\paragraph{IV. Radial instance.}
For \(\alpha>B\epsilon\), rescaling the rotating construction of Instance~II and, as
in~\cite{li2026smoothing}, initializing on its rotating eigenmode removes the transient term for every
\(\beta\in[0,1]\), keeping the iterates on a circle with non-decaying KKT
residuals; the explicit instance and proof are given in
the
supplementary material.

\begin{lemma}
\label{lem:main-radial}
There exist constants \(B\ge1\), \(\epsilon_{\mathrm{r}}\in(0,\tfrac12]\), and
\(c_{\mathrm{r}}>0\) such that, for every \(\epsilon\in(0,\epsilon_{\mathrm{r}})\),
\(\alpha\in(B\epsilon,1]\), and \(\beta\in[0,1]\), there exist an
instance \(I\in\mathcal I(\mathbf C,\Delta)\) and an initialization
\(\bu^0\in\mathcal S_0(\mathbf C)\) such that
\(\mathrm{LRALM}_{\theta(\epsilon),\alpha,\beta}\) initialized at \(\bu^0\)
satisfies \(R_k\ge c_{\mathrm{r}}\alpha>\epsilon\) for all \(k\ge0\); hence
\(T_\epsilon(\theta(\epsilon),\alpha,\beta;I,\bu^0)=+\infty\).
\end{lemma}

Lemmas~\ref{lem:main-plateau} and~\ref{lem:main-radial} together yield
Theorem~\ref{thm:lower-bound-fixed-rho}, whose complete proof appears in
the supplementary material.

\section{Acceleration Beyond the Family with Fast LRALM}
\label{sec:acceleration}

By Theorem~\ref{thm:lower-bound-fixed-rho} and
Proposition~\ref{prop:upper-beta0}, \(\Theta(\epsilon^{-3})\) is the best
order that \(\calC_\theta\) attains with a constant \(\rho\) and
\(s=c_s/\epsilon\), so a faster rate requires leaving the family. To address
\hyperlink{Q2}{\textbf{(Q2)}}, we construct Fast LRALM, a constant-penalty
single-loop method that performs inexact Nesterov accelerated gradient ascent on
the regularized dual objective, which is smooth and strongly concave under
a constant penalty and primal regularization. Throughout
this section, we set \(\bmu=\zero\) and write
\(\mathcal P_s(\bx,\blmbd;\bz):=\mathcal P_s(\bx,\blmbd;\bz,\zero)\). In
Table~\ref{tab:acc-notation}, we summarize the value functions and solution
mappings, which are well defined under \(\gamma L_{\calD,\bx}<1/2\) by
Proposition~\ref{prop:curvature-smoothness}.

\begin{table}[t]
\centering
\small
\caption{Value functions and solution mappings for Fast LRALM.}
\label{tab:acc-notation}
\renewcommand{\arraystretch}{1.2}
\setlength{\tabcolsep}{4pt}
\begin{tabular}{@{}lll@{}}
\toprule
Notation & Definition & Notes \\
\midrule
\(\mathcal P_s(\bx,\blmbd;\bz)\)
  & see~\eqref{eq:potential} with \(\bmu=\zero\)
  & regularized augmented Lagrangian \\
\(\phi_s(\bx)\)
  & \(\max_{\blmbd}\bigl(\calD_\rho(\bx,\blmbd)-\frac{1}{2s}\|\blmbd\|^2\bigr)\)
  & primal function \\
\(d_s(\blmbd;\bz)\) & \(\min_{\bx}\mathcal P_s(\bx,\blmbd;\bz)\) & dual function \\
\(p_s(\bz)\) & \(\max_{\blmbd}\min_{\bx}\mathcal P_s(\bx,\blmbd;\bz)\)
  & proximal function; optimal value \\
\(\bx_s(\blmbd;\bz)\) & \(\argmin_{\bx}\mathcal P_s(\bx,\blmbd;\bz)\) & primal minimizer \\
\(\blmbd_s^\star(\bz)\) & \(\argmax_{\blmbd}d_s(\blmbd;\bz)\) & dual optimal solution \\
\(\bx_s^\star(\bz)\) & \(\bx_s(\blmbd_s^\star(\bz);\bz)\) & \(=\prox_{\gamma\phi_s}(\bz)\) \\
\bottomrule
\end{tabular}
\end{table}

We target the three-level coupled problem
\begin{equation}
    \min_{\bz\in\RR^n}\left\{p_s(\bz)=\max_{\blmbd\in\RR^m}\min_{\bx\in\RR^n}
    \mathcal P_s(\bx,\blmbd;\bz)\right\}.
    \label{eq:fast-lralm-nested}
\end{equation}
Proposition~\ref{prop:curvature-smoothness}(b) and the definition of \(\phi_s\)
give
\begin{equation}
    p_s(\bz)
    =\min_{\bx\in\RR^n}\max_{\blmbd\in\RR^m}\mathcal P_s(\bx,\blmbd;\bz)
    =\min_{\bx\in\RR^n}\Bigl\{\phi_s(\bx)+\frac{1}{2\gamma}\|\bx-\bz\|^2\Bigr\},
    \label{eq:minimax-equality}
\end{equation}
which implies that \(p_s\) is the Moreau envelope of the
\(L_{\calD,\bx}\)-weakly convex function \(\phi_s\). Under
\(\gamma L_{\calD,\bx}<1/2\), \(p_s\) is smooth with
\(\nabla p_s(\bz)=(\bz-\bx_s^\star(\bz))/\gamma\), where
\(\bx_s^\star(\bz)=\bx_s(\blmbd_s^\star(\bz);\bz)\). We therefore seek an approximate stationary point
of \(p_s\), i.e., a center \(\bz\) with \(\|\bz-\bx_s^\star(\bz)\|\) small.

Algorithm~\ref{alg:b2-prime} summarizes Fast LRALM, which
solves~\eqref{eq:fast-lralm-nested} in a single loop. At each iterate \(\bz^k\),
it takes one inexact gradient step on \(p_s\), together with one dual step and
one primal step that track the unique saddle point
\(\bigl(\bx_s^\star(\bz^k),\blmbd_s^\star(\bz^k)\bigr)\) of the strongly
convex-strongly concave inner problem in~\eqref{eq:fast-lralm-nested}. We
describe these steps from the outer level to the inner one.

\begin{algorithm}[t]
    \caption{Fast LRALM}
    \label{alg:b2-prime}
    \small
    \begin{algorithmic}[1]
        \STATE \textbf{Input:} Parameters \(\rho,s,\gamma,\eta,\tau>0\),
        \(\alpha\in(0,1]\); \(\kappa:=s/\tau\),
        \(\omega:=\frac{\sqrt{\kappa}-1}{\sqrt{\kappa}+1}\); iteration limit \(K\).
        \STATE \textbf{Initialize:} \((\bx^0,\blmbd^0,\bz^0,\zero)\in\mathcal S_0(\mathbf C)\)
        and \(\blmbd^{-1}=\blmbd^0\).
        \FOR{\(k=0,1,\ldots,K-1\)}
            \STATE \(\bar{\blmbd}^k=\blmbd^k+\omega(\blmbd^k-\blmbd^{k-1})\).
            \STATE \(\bx^{k+1}=\bx^k-\eta\Bigl(\nabla_{\bx}\calL_\rho\bigl(\bx^k,\by(\bx^k,\bar{\blmbd}^k),\bar{\blmbd}^k\bigr)+\tfrac{1}{\gamma}(\bx^k-\bz^k)\Bigr)\).
            \STATE \(\blmbd^{k+1}=\bar{\blmbd}^k+\tau\Bigl(\calA(\bx^{k+1})-\by(\bx^{k+1},\bar{\blmbd}^k)-\tfrac{1}{s}\bar{\blmbd}^k\Bigr)\).
            \STATE \(\bz^{k+1}=\bz^k+\alpha(\bx^{k+1}-\bz^k)\).
        \ENDFOR
        \STATE \textbf{Output:}
        \(\bigl\{\bigl(\bx^{k+1},\,\by(\bx^{k+1},\bar{\blmbd}^k),\,\bq(\bx^{k+1},\bar{\blmbd}^k)\bigr)\bigr\}_{k=0}^{K-1}\).
    \end{algorithmic}
\end{algorithm}
\paragraph{Inexact gradient descent on \(p_s\).}
Since \(p_s\) is smooth, a natural approach is gradient descent on \(p_s\).
However, the gradient \(\nabla p_s(\bz^k)=(\bz^k-\bx_s^\star(\bz^k))/\gamma\)
involves the minimizer \(\bx_s^\star(\bz^k)\) of~\eqref{eq:minimax-equality},
whose exact computation requires solving a minimax subproblem and leads to
double-loop methods~\cite{kong2021accelerated,lin2020near}. Following
the single-loop proximal-center framework~\cite{zhang2020smoothed,li2026smoothing,zhang2026matching},
we replace \(\bx_s^\star(\bz^k)\) with the primal iterate \(\bx^{k+1}\), which
approximates \(\bx_s^\star(\bz^k)\) through the dual and primal steps below, and take
one inexact gradient step with stepsize \(\alpha\gamma\), giving
\begin{equation}
    \bz^{k+1}=\bz^k-\alpha\gamma\,\frac{\bz^k-\bx^{k+1}}{\gamma}
    =\bz^k+\alpha(\bx^{k+1}-\bz^k).
    \label{eq:fast-lralm-center-step}
\end{equation}
\paragraph{Nesterov's accelerated gradient ascent on \(d_s\).}
Since \(\bx_s^\star(\bz^k)=\bx_s(\blmbd_s^\star(\bz^k);\bz^k)\), the exact center
step requires solving the dual maximization problem \(\max_{\blmbd}d_s(\blmbd;\bz^k)\). For every \(\bz\in\RR^n\),
\(d_s(\cdot;\bz)\) is \(s^{-1}\)-strongly concave and smooth with a modulus
independent of \(s\), as established in the supplementary material, so its condition number
grows with \(s\). We therefore take one step of Nesterov's
accelerated gradient ascent~\cite{nesterov2018lectures} with parameters
\(\kappa\) and \(\omega\) as in Algorithm~\ref{alg:b2-prime}. By Danskin's theorem,
\(\nabla d_s(\bar{\blmbd}^k;\bz^k)
=\nabla_{\blmbd}\mathcal P_s(\bx_s(\bar{\blmbd}^k;\bz^k),\bar{\blmbd}^k;\bz^k)\);
replacing the exact minimizer \(\bx_s(\bar{\blmbd}^k;\bz^k)\) with the primal
iterate \(\bx^{k+1}\) gives
\begin{equation}
    \blmbd^{k+1}
    =\bar{\blmbd}^k+\tau\nabla_{\blmbd}\mathcal P_s(\bx^{k+1},\bar{\blmbd}^k;\bz^k)
    =\bar{\blmbd}^k+\tau\Bigl(\calA(\bx^{k+1})-\by(\bx^{k+1},\bar{\blmbd}^k)
    -\frac1s\bar{\blmbd}^k\Bigr).
    \label{eq:fast-lralm-dual-step}
\end{equation}
Compared with Algorithm~\ref{alg:lralm}, the extrapolation \(\omega>0\) is the
sole departure from \(\calC_\theta\).

\paragraph{Primal gradient descent on \(\mathcal P_s\).}
The exact dual gradient requires the minimizer
\(\bx_s(\bar{\blmbd}^k;\bz^k)=\argmin_{\bx\in\RR^n}\mathcal P_s(\bx,\bar{\blmbd}^k;\bz^k)\).
Proposition~\ref{prop:curvature-smoothness}(a) ensures that the condition number
\(L_x/\mu_\gamma\) of \(\mathcal P_s(\cdot,\bar{\blmbd}^k;\bz^k)\) is
independent of \(\bar{\blmbd}^k\) and \(s\); hence we take a
single gradient descent step from \(\bx^k\) to track
\(\bx_s(\bar{\blmbd}^k;\bz^k)\), and ultimately \(\bx_s^\star(\bz^k)\), giving
\begin{equation}
    \bx^{k+1}
    =\bx^k-\eta\nabla_{\bx}\mathcal P_s(\bx^k,\bar{\blmbd}^k;\bz^k)
    =\bx^k-\eta\Bigl(\nabla_{\bx}\calL_\rho\bigl(\bx^k,\by(\bx^k,\bar{\blmbd}^k),
    \bar{\blmbd}^k\bigr)+\frac1\gamma(\bx^k-\bz^k)\Bigr),
    \label{eq:fast-lralm-primal-step}
\end{equation}
where the second equality follows from~\eqref{eq:reduced-al-gradient} and
\(\nabla_{\bx}\calL_\rho(\bx,\by(\bx,\blmbd),\blmbd)
=\nabla f(\bx)+J_{\calA}(\bx)^\top\bq(\bx,\blmbd)\).

Since each step uses the estimate produced by the level below it, the updates
are executed from the innermost level outward by forming \(\bar{\blmbd}^k\) and
applying~\eqref{eq:fast-lralm-primal-step}, \eqref{eq:fast-lralm-dual-step},
and~\eqref{eq:fast-lralm-center-step} in turn. Substituting the gradients
in~\eqref{eq:reduced-al-gradient} yields Algorithm~\ref{alg:b2-prime}.

\begin{remark}
\begin{enumerate}[label=(\roman*)]
\item The dual step~\eqref{eq:fast-lralm-dual-step} can be written as the damped
multiplier update
\[
    \blmbd^{k+1}=(1-\tau/s)\bar{\blmbd}^k
    +\tau\bigl(\calA(\bx^{k+1})-\by(\bx^{k+1},\bar{\blmbd}^k)\bigr),
\]
which resembles the multiplier update of single-loop perturbed ALMs such as
GDPA~\cite{luSingleLoopGradientDescent2022} and
LIPAL~\cite{elbourkhissi2025perturbed}, except that it is evaluated at the
extrapolated multiplier \(\bar{\blmbd}^k\).
\item Without the primal
regularization term \(\frac{1}{2\gamma}\|\bx-\bz\|^2\),
\(\calD_\rho(\cdot,\blmbd)\) is nonconvex and need not induce a smooth dual
objective for Nesterov extrapolation; including this term makes
\(d_s(\cdot;\bz^k)\) smooth and strongly concave. Inexact
ALMs~\cite{sahin2019inexact,li2021rate} solve outer augmented Lagrangian
subproblems inside an inner loop until a stopping criterion is met. In contrast,
Fast LRALM moves \(\bz^k\) at every iteration with stepsize
\(\alpha=\Theta(\kappa^{-1/2})\), so \(d_s(\cdot;\bz^k)\) changes little
between iterations and a single primal-dual step tracks the saddle point without
inner loops.
\end{enumerate}
\end{remark}

\subsection{KKT residual bounds}
\label{subsec:fast-lralm-residual}
We now analyze the iteration complexity of Fast LRALM. Throughout the
remainder of this section, let \(\mathbf C\in\RR^5_{\ge0}\), \(\Delta>0\),
\(I=(f,\calA,h)\in\mathcal I(\mathbf C,\Delta)\), and
\(\bu^0=(\bx^0,\blmbd^0,\bz^0,\zero)\in\mathcal S_0(\mathbf C)\). Let
\(\{(\bx^k,\blmbd^k,\bz^k)\}_{k\ge0}\) be the sequence generated by
Algorithm~\ref{alg:b2-prime} from \(\bu^0\), and write its output as
\((\bx^{k+1},\by^{k+1},\bq^k)\), where
\begin{equation}
    \by^{k+1}:=\by(\bx^{k+1},\bar{\blmbd}^k),
    \qquad \bq^k:=\bq(\bx^{k+1},\bar{\blmbd}^k).
    \label{eq:b2-prime-outputs}
\end{equation}
The following lemma bounds the KKT residual at \((\bx^{k+1},\by^{k+1},\bq^k)\).
\begin{lemma}
\label{lem:b2-prime-residual-decomposition}
Suppose that Assumptions~\ref{ass:smoothness} and~\ref{ass:h} hold. Then,
for every \(k\ge0\), \(\bq^k\in\partial h(\by^{k+1})\), and the two components
of \(\mathcal R_{\mathrm{KKT}}(\bx^{k+1},\by^{k+1},\bq^k)\)
in~\eqref{eq:kkt-residual} satisfy
\[
\begin{aligned}
    \bigl\|\nabla f(\bx^{k+1})+J_{\calA}(\bx^{k+1})^\top\bq^k\bigr\|
    &\le\left(\left|\frac1\eta-\frac1\gamma\right|+L_{\calD,\bx}\right)\|\bx^{k+1}-\bx^k\|
    +\frac1\gamma\|\bx^{k+1}-\bz^k\|,\\
    \|\calA(\bx^{k+1})-\by^{k+1}\|
    &\le\frac1\tau\|\blmbd^{k+1}-\bar{\blmbd}^k\|
    +\frac1s\|\bar{\blmbd}^k\|.
\end{aligned}
\]
\end{lemma}

\begin{proof}
Lemma~\ref{lem:reduced-al}(a) gives \(\bq^k\in\partial h(\by^{k+1})\) and
\(\nabla_{\bx}\calD_\rho(\bx^{k+1},\bar{\blmbd}^k)
=\nabla f(\bx^{k+1})+J_{\calA}(\bx^{k+1})^\top\bq^k\). For the stationarity
component, the primal step of Algorithm~\ref{alg:b2-prime}
and~\eqref{eq:potential} yield
\[
    \frac{\bx^k-\bx^{k+1}}{\eta}
    =\nabla_{\bx}\calD_\rho(\bx^k,\bar{\blmbd}^k)+\frac{1}{\gamma}(\bx^k-\bz^k).
\]
Rearranging these two identities gives
\[
\begin{aligned}
    \nabla f(\bx^{k+1})+J_{\calA}(\bx^{k+1})^\top\bq^k
    &=\bigl[\nabla_{\bx}\calD_\rho(\bx^{k+1},\bar{\blmbd}^k)
    -\nabla_{\bx}\calD_\rho(\bx^k,\bar{\blmbd}^k)\bigr]\\
    &\quad+\Bigl(\frac{1}{\eta}-\frac{1}{\gamma}\Bigr)(\bx^k-\bx^{k+1})
    -\frac{1}{\gamma}(\bx^{k+1}-\bz^k).
\end{aligned}
\]
As shown in Section~\ref{subsec:regularized-al},
\(\nabla_{\bx}\calD_\rho(\cdot,\bar{\blmbd}^k)\) is
\(L_{\calD,\bx}\)-Lipschitz, so the triangle inequality gives the first
bound. For the feasibility component, the dual step of
Algorithm~\ref{alg:b2-prime} and~\eqref{eq:b2-prime-outputs} give
\[
    \frac{\blmbd^{k+1}-\bar{\blmbd}^k}{\tau}
    =\calA(\bx^{k+1})-\by^{k+1}-\frac{1}{s}\bar{\blmbd}^k,
\]
from which the second bound follows by the triangle inequality, completing the proof.
\end{proof}

Since Fast LRALM takes unprojected dual steps with extrapolation
coefficient \(1+\omega>1\), multiplier boundedness does not follow from
projection or convex combination and must be established separately. We next show that the multipliers remain bounded, adapting the quadratic Lyapunov function of~\cite{khlebnikov2011optimization} to the dual recursion of Fast LRALM.
\begin{lemma}
\label{lem:b2-prime-multiplier-bounds}
Suppose that Assumption~\ref{ass:h} holds and
\(\tau L_\lambda\le\frac12\). Then, for every \(k\ge0\),
\[
    \max\left\{\|\blmbd^k\|,\|\bar{\blmbd}^k\|\right\}
    \le6\sqrt{\frac{\rho}{\tau}}\,L_h.
\]
\end{lemma}

\begin{proof}
Write \(t:=\tau L_\lambda\). Since~\eqref{eq:reduced-al-gradient} gives
\(\calA(\bx^{k+1})-\by^{k+1}=(\bq^k-\bar{\blmbd}^k)/\rho\), the dual
step of Algorithm~\ref{alg:b2-prime} becomes
\[
    \blmbd^{k+1}=(1-t)\bar{\blmbd}^k+\frac{\tau}{\rho}\bq^k,
    \qquad
    \bar{\blmbd}^k=(1+\omega)\blmbd^k-\omega\blmbd^{k-1},
    \qquad \|\bq^k\|\le L_h.
\]
With \(a:=(1-t)(1+\omega)\) and \(b:=(1-t)\omega\), this is the two-step
recursion
\begin{equation}
    \blmbd^{k+1}=a\blmbd^k-b\blmbd^{k-1}+\frac{\tau}{\rho}\bq^k .
    \label{eq:b2-prime-two-step}
\end{equation}
Because \(\tau/s\le t\le\frac12\), we have \(\omega\in(0,1)\), so
\(0<b\le1-t\). Consider, for \(k\ge0\), the quadratic function
\[
    V_k:=\|\blmbd^k\|^2
    -a\langle\blmbd^k,\blmbd^{k-1}\rangle
    +b\|\blmbd^{k-1}\|^2 .
\]
Completing the square gives
\[
    V_k=\Bigl(1-\frac{a^2}{4b}\Bigr)\|\blmbd^k\|^2
    +\frac{1}{4b}\|a\blmbd^k-2b\blmbd^{k-1}\|^2 .
\]
The definition of \(\omega\) gives
\((1-\omega)^2=\frac{\tau}{s}(1+\omega)^2\), which, combined with
\(t-\frac{\tau}{s}=\frac{\tau}{\rho}\), leads to
\[
    1-\frac{a^2}{4b}
    =\frac{t(1+\omega)^2-(1-\omega)^2}{4\omega}
    =\frac{(1+\omega)^2}{4\omega}\cdot\frac{\tau}{\rho}
    \ge\frac{\tau}{\rho},
\]
where the inequality follows from \((1+\omega)^2\ge4\omega\). Hence both terms
in the completed square are nonnegative, and
\begin{equation}
    V_k\ge\frac{\tau}{\rho}\|\blmbd^k\|^2,
    \qquad
    V_k\ge\frac{1}{4b}\|a\blmbd^k-2b\blmbd^{k-1}\|^2 .
    \label{eq:b2-prime-V-squares}
\end{equation}
To propagate \(V_k\) along the iterations, we substitute~\eqref{eq:b2-prime-two-step}
into the definition of \(V_{k+1}\) and obtain
\begin{align*}
    V_{k+1}
    &=\langle\blmbd^{k+1},\blmbd^{k+1}-a\blmbd^k\rangle+b\|\blmbd^k\|^2\\
    &=\Bigl\langle a\blmbd^k-b\blmbd^{k-1}+\frac{\tau}{\rho}\bq^k,\,
    -b\blmbd^{k-1}+\frac{\tau}{\rho}\bq^k\Bigr\rangle+b\|\blmbd^k\|^2\\
    &=bV_k+\frac{\tau}{\rho}\langle\bq^k,a\blmbd^k-2b\blmbd^{k-1}\rangle
    +\frac{\tau^2}{\rho^2}\|\bq^k\|^2\\
    &\le bV_k+2\sqrt{bV_k}\,\frac{\tau}{\rho}\|\bq^k\|
    +\frac{\tau^2}{\rho^2}\|\bq^k\|^2
    =\Bigl(\sqrt{bV_k}+\frac{\tau}{\rho}\|\bq^k\|\Bigr)^2,
\end{align*}
where the inequality follows from the Cauchy--Schwarz inequality and the second
inequality in~\eqref{eq:b2-prime-V-squares}. Taking square roots, and using \(\|\bq^k\|\le L_h\),
\(\sqrt{b}\le\sqrt{1-t}\le1-t/2\), and \(\tau/\rho\le t\), we arrive at
\[
    \sqrt{V_{k+1}}
    \le\sqrt{b}\,\sqrt{V_k}+\frac{\tau}{\rho}L_h
    \le\Bigl(1-\frac t2\Bigr)\sqrt{V_k}+\frac t2\cdot2L_h .
\]
Since \(\blmbd^{-1}=\blmbd^0\), \(\|\blmbd^0\|\le L_h\), and \(t\le\frac12\),
we have
\[
    V_0=(1-a+b)\|\blmbd^0\|^2=t\|\blmbd^0\|^2\le\frac12L_h^2,
\]
so \(\sqrt{V_0}\le 2L_h\). Induction on the preceding recursion then gives
\(\sqrt{V_k}\le 2L_h\) for all \(k\ge0\). Combining this bound with the first
inequality in~\eqref{eq:b2-prime-V-squares} and using \(\omega\in(0,1)\) yield
\[
    \|\blmbd^k\|\le\sqrt{\frac{\rho}{\tau}}\,\sqrt{V_k}
    \le 2\sqrt{\frac{\rho}{\tau}}\,L_h,
    \qquad
    \|\bar{\blmbd}^k\|\le(1+\omega)\|\blmbd^k\|+\omega\|\blmbd^{k-1}\|
    \le 6\sqrt{\frac{\rho}{\tau}}\,L_h,
\]
which completes the proof.
\end{proof}

\subsection{Analysis framework}
\label{subsec:fast-lralm-lyapunov}
We present the key ideas of the analysis in three steps.
\begin{stepbox}
\textbf{Step 1:} Construct a Lyapunov sequence \(\{\mathcal W_k\}_{k\ge0}\)
whose one-step decrease controls the increments of the primal step in \(\bx\),
the dual step in \(\blmbd\), and the center step in \(\bz\).
\end{stepbox}

Building on the potential function
in~\cite{zhang2020smoothed,zheng2023universal,li2025nonsmooth,li2026smoothing}, we
define the Lyapunov sequence tied to Algorithm~\ref{alg:b2-prime} as
\[
    \mathcal W_k:=\underbrace{p_s(\bz^k)}_{\text{center}}
    +\frac12\mathcal E_k(\bz^k),\qquad k\ge0,
\]
where
\[
    \mathcal E_k(\bz)
    :=\underbrace{p_s(\bz)-d_s(\blmbd^k;\bz)
    +\frac{1}{2s}\|\bw^k-\blmbd_s^\star(\bz)\|^2}_{\text{dual}}
    +\underbrace{\frac{\mu_\gamma}{4}\|\bx^k-\bx_s(\bar{\blmbd}^{k-1};\bz)\|^2
    +\frac{G_{\calA}^2}{4\mu_\gamma}\|\blmbd^k-\bar{\blmbd}^{k-1}\|^2}_{\text{primal}},
\]
with \(\bw^k:=\blmbd^{k-1}+\sqrt{\kappa}(\blmbd^k-\blmbd^{k-1})\),
\(\bar{\blmbd}^{-1}:=\blmbd^0\), and \(\mathcal E_k:=\mathcal E_k(\bz^k)\).
The above Lyapunov sequence differs from this potential function in the dual and
primal parts. The dual part adds
the distance term of the estimate-sequence Lyapunov
function~\cite{nesterov2018lectures,wilson2021lyapunov}, since under Nesterov
extrapolation the dual gap alone need not decrease. The primal part replaces
the primal gap \(\mathcal P_s(\bx^k,\blmbd^k;\bz)-d_s(\blmbd^k;\bz)\) with two
quadratic terms, which control the error of evaluating
\(\nabla_{\blmbd}\mathcal P_s(\bx^{k+1},\bar{\blmbd}^k;\bz^k)\) in place of
\(\nabla d_s(\bar{\blmbd}^k;\bz^k)\).
Following the order of the updates in Algorithm~\ref{alg:b2-prime}, we
decompose the one-step change of \(\mathcal W_k\) as
\begin{equation}
    \mathcal W_{k+1}-\mathcal W_k
    =\frac12\underbrace{\bigl[\mathcal E_{k+1}(\bz^k)-\mathcal E_k(\bz^k)\bigr]}_{\text{primal-dual contraction}}
    +\underbrace{\bigl[p_s(\bz^{k+1})-p_s(\bz^k)\bigr]}_{\text{proximal descent}}
    +\frac12\underbrace{\bigl[\mathcal E_{k+1}(\bz^{k+1})-\mathcal E_{k+1}(\bz^k)\bigr]}_{\text{proximal-point shift}}.
    \label{eq:b2-prime-decomposition}
\end{equation}

We bound the three terms in~\eqref{eq:b2-prime-decomposition} in turn,
starting with the first two.
\begin{stepbox}
\textbf{Step 2:} Bound the primal-dual contraction and proximal-descent
terms in~\eqref{eq:b2-prime-decomposition} via
Lemmas~\ref{lem:b2-prime-fixed-center-contraction}
and~\ref{lem:b2-prime-z-step-p}.
\end{stepbox}

With \(\bz^k\) held fixed, Lemma~\ref{lem:b2-prime-fixed-center-contraction}
shows that one primal-dual step with dual Nesterov extrapolation contracts
\(\mathcal E_k\) at the accelerated rate \(1-\kappa^{-1/2}\) while producing
negative squared primal and dual increments together with a negative
primal-subproblem term \(-\Theta(\kappa^{-1/2})\|\bx^{k+1}-\bx_s^\star(\bz^k)\|^2\).
For the proximal-descent term, the center step~\eqref{eq:fast-lralm-center-step}
is an inexact gradient step on \(p_s\) that replaces \(\bx_s^\star(\bz^k)\)
with \(\bx^{k+1}\); by Lemma~\ref{lem:b2-prime-z-step-p}, it therefore yields a
negative center term \(-\Theta(\alpha)\|\bx^{k+1}-\bz^k\|^2\) at the cost of a
positive primal-subproblem error \(\Theta(\alpha)\|\bx^{k+1}-\bx_s^\star(\bz^k)\|^2\).

We now turn to the remaining proximal-shift term
in~\eqref{eq:b2-prime-decomposition}. Since the center moves at every
iteration without an inner loop, the resulting increase in \(\mathcal E_{k+1}\)
is absorbed by the contraction and negative center term from Step~2, as
in~\cite{zhang2026matching}.

\begin{stepbox}
\textbf{Step 3:} Bound the proximal-shift term and choose
\(\alpha=\Theta(\kappa^{-1/2})\) to balance the outer center step against the
inner contraction rate via Lemmas~\ref{lem:b2-prime-center-error}
and~\ref{lem:b2-prime-potential-descent}.
\end{stepbox}

Specifically, Lemma~\ref{lem:b2-prime-center-error} shows that moving \(\bz^k\)
to \(\bz^{k+1}\) shifts the inner saddle point and increases \(\mathcal E_{k+1}\)
by an amount controlled by \(\alpha\|\bx^{k+1}-\bz^k\|\), with constants
independent of \(s\). If \(\alpha\) is too large, this drift outpaces the
inner contraction and prevents descent. Matching the center stepsize to the
accelerated inner contraction rate, \(\alpha=\Theta(\kappa^{-1/2})\),
as in~\cite{zhang2026matching}, ensures
that each iterate warm-starts the next subproblem; the contraction margin and
negative terms from Step~2 then absorb both the proximal-shift increase and the
primal-subproblem error, yielding the sufficient decrease of
\(\{\mathcal W_k\}_{k\ge0}\) in Lemma~\ref{lem:b2-prime-potential-descent}. At
\(s=\Theta(\epsilon^{-1})\), this yields \(\alpha=\Theta(\epsilon^{1/2})\), larger
than the \(\Theta(\epsilon)\) threshold in \(\calC_\theta\).

\subsection{Iteration complexity analysis}
\label{subsec:fast-lralm-analysis}
Throughout this subsection, we impose the parameter conditions
\begin{equation}
    \eta=\gamma,\quad \gamma L_{\calD,\bx}\le\frac{1}{10},\quad
    0<\tau\le\frac14\min\{\rho,s\},\quad 14\gamma\tau G_{\calA}^2\le1,\quad
    \alpha=\frac{c_\alpha}{\sqrt{\kappa}},\quad c_\alpha\in\Bigl(0,\frac16\Bigr],
    \label{eq:b2-prime-stepsizes}
\end{equation}
and use the following properties of the value functions in
Table~\ref{tab:acc-notation}. The Lipschitz bounds on \(\bx_s(\cdot;\cdot)\)
and \(\bx_s^\star(\cdot)\) parallel those in~\cite{zhang2020smoothed}, and
the bound on \(\blmbd_s^\star(\cdot)\) additionally uses the
\(s^{-1}\)-strong concavity of \(\mathcal P_s\) in \(\blmbd\); the proofs
are given in the supplementary material.

\begin{lemma}
\label{lem:b2-prime-fixed-center-values}
Suppose that Assumptions~\ref{ass:smoothness} and~\ref{ass:h} hold. Then, for
every \(\bz\in\RR^n\), the following statements hold.
\begin{enumerate}[label=(\roman*)]
\item For all \(\blmbd,\blmbd'\in\RR^m\),
\[
    \|\bx_s(\blmbd;\bz)-\bx_s(\blmbd';\bz)\|
    \le\frac{G_{\calA}}{\mu_\gamma}\|\blmbd-\blmbd'\| .
\]
\item For all \(\bx\in\RR^n\) and \(\blmbd\in\RR^m\),
\[
    \frac{\mu_\gamma}{2}\|\bx-\bx_s(\blmbd;\bz)\|^2
    \le\mathcal P_s(\bx,\blmbd;\bz)-d_s(\blmbd;\bz)
    \le\frac{L_x}{2}\|\bx-\bx_s(\blmbd;\bz)\|^2 .
\]
\item For every \(\blmbd\in\RR^m\),
\[
    p_s(\bz)-d_s(\blmbd;\bz)
    \ge\frac{1}{2s}\|\blmbd-\blmbd_s^\star(\bz)\|^2
    +\frac{\mu_\gamma}{2}\|\bx_s(\blmbd;\bz)-\bx_s^\star(\bz)\|^2 .
\]
\end{enumerate}
\end{lemma}

\begin{lemma}
\label{lem:b2-prime-center-dependence}
Suppose that Assumptions~\ref{ass:smoothness} and~\ref{ass:h} hold. Let
\(\bz,\bz'\in\RR^n\) and \(\blmbd\in\RR^m\). Then
\begin{align*}
    p_s(\bz')-p_s(\bz)
    &\le\frac1\gamma\bigl\langle\bz-\bx_s^\star(\bz),\bz'-\bz\bigr\rangle
    +\frac{1}{2\gamma}\|\bz'-\bz\|^2,\\
    d_s(\blmbd;\bz')-d_s(\blmbd;\bz)
    &\ge\frac1\gamma\bigl\langle\bz-\bx_s(\blmbd;\bz'),\bz'-\bz\bigr\rangle
    +\frac{1}{2\gamma}\|\bz'-\bz\|^2 ,
\end{align*}
and
\begin{gather*}
    \max\bigl\{\|\bx_s(\blmbd;\bz')-\bx_s(\blmbd;\bz)\|,\,
    \|\bx_s^\star(\bz')-\bx_s^\star(\bz)\|\bigr\}
    \le\frac{\|\bz'-\bz\|}{\gamma\mu_\gamma},\\
    \|\blmbd_s^\star(\bz')-\blmbd_s^\star(\bz)\|
    \le\frac{\sqrt{s}\,\|\bz'-\bz\|}{2\gamma\sqrt{\mu_\gamma}} .
\end{gather*}
\end{lemma}

We begin with Step~2 and bound
the primal-dual contraction term in~\eqref{eq:b2-prime-decomposition}.

\begin{lemma}
\label{lem:b2-prime-fixed-center-contraction}
Suppose that Assumptions~\ref{ass:smoothness} and~\ref{ass:h}
and~\eqref{eq:b2-prime-stepsizes} hold. Then, for every \(k\ge0\),
\[
    \mathcal E_{k+1}(\bz^k)
    \le\Bigl(1-\frac{1}{\sqrt{\kappa}}\Bigr)\mathcal E_k
    -\frac{1}{64\gamma}\|\bx^{k+1}-\bx^k\|^2
    -\frac{1}{8\tau}\|\blmbd^{k+1}-\bar{\blmbd}^k\|^2
    -\frac{\mu_\gamma}{2\sqrt{\kappa}}\|\bx^{k+1}-\bx_s^\star(\bz^k)\|^2 .
\]
\end{lemma}

\begin{proof}
\emph{Step 1 (Dual accelerated inequality).} The definitions of \(\bw^k\), \(\omega\), and \(\bar{\blmbd}^k\) give the identity
\[
    \bw^{k+1}-\sqrt{\kappa}(\blmbd^{k+1}-\bar{\blmbd}^k)
    =\Bigl(1-\frac{1}{\sqrt{\kappa}}\Bigr)\bw^k
    +\frac{1}{\sqrt{\kappa}}\bar{\blmbd}^k
    =\bar{\blmbd}^k+(\sqrt{\kappa}-1)(\bar{\blmbd}^k-\blmbd^k).
\]
Expanding \(\|\bw^{k+1}-\blmbd_s^\star(\bz^k)\|^2\) with this identity, and
invoking the convexity of \(\|\cdot-\blmbd_s^\star(\bz^k)\|^2\) together with
\(\kappa/s=1/\tau\), we get
\begin{equation}
\begin{aligned}
    \frac{1}{2s}\|\bw^{k+1}-\blmbd_s^\star(\bz^k)\|^2
    &\le\Bigl(1-\frac{1}{\sqrt{\kappa}}\Bigr)
    \frac{1}{2s}\|\bw^k-\blmbd_s^\star(\bz^k)\|^2
    +\frac{1}{2s\sqrt{\kappa}}\|\bar{\blmbd}^k-\blmbd_s^\star(\bz^k)\|^2\\
    &\quad+\frac{1}{2\tau}\|\blmbd^{k+1}-\bar{\blmbd}^k\|^2
    +\frac{1}{\tau\sqrt{\kappa}}
    \bigl\langle\blmbd^{k+1}-\bar{\blmbd}^k,\,
    (\sqrt{\kappa}-1)(\bar{\blmbd}^k-\blmbd^k)
    +\bar{\blmbd}^k-\blmbd_s^\star(\bz^k)\bigr\rangle .
\end{aligned}
\label{eq:b2-prime-w-expand}
\end{equation}
The dual step of Algorithm~\ref{alg:b2-prime} reads
\(\blmbd^{k+1}-\bar{\blmbd}^k
=\tau\nabla_{\blmbd}\mathcal P_s(\bx^{k+1},\bar{\blmbd}^k;\bz^k)\). As
\(\mathcal P_s(\bx^{k+1},\cdot;\bz^k)\) is \(s^{-1}\)-strongly concave and
\(L_\lambda\)-smooth, the inner products on the right-hand side admit the bounds
\begin{align*}
    \frac1\tau\bigl\langle\blmbd^{k+1}-\bar{\blmbd}^k,
    \bar{\blmbd}^k-\blmbd^k\bigr\rangle
    &\le\mathcal P_s(\bx^{k+1},\bar{\blmbd}^k;\bz^k)
    -\mathcal P_s(\bx^{k+1},\blmbd^k;\bz^k),\\
    \frac1\tau\bigl\langle\blmbd^{k+1}-\bar{\blmbd}^k,
    \bar{\blmbd}^k-\blmbd_s^\star(\bz^k)\bigr\rangle
    &\le\mathcal P_s(\bx^{k+1},\bar{\blmbd}^k;\bz^k)
    -\mathcal P_s(\bx^{k+1},\blmbd_s^\star(\bz^k);\bz^k)\\
    &\quad-\frac{1}{2s}\|\bar{\blmbd}^k-\blmbd_s^\star(\bz^k)\|^2,\\
    \mathcal P_s(\bx^{k+1},\bar{\blmbd}^k;\bz^k)
    &\le\mathcal P_s(\bx^{k+1},\blmbd^{k+1};\bz^k)
    -\frac{2-\tau L_\lambda}{2\tau}\|\blmbd^{k+1}-\bar{\blmbd}^k\|^2 .
\end{align*}
Inserting these three inequalities into~\eqref{eq:b2-prime-w-expand} and
adding \(p_s(\bz^k)-d_s(\blmbd^{k+1};\bz^k)
-(1-\kappa^{-1/2})[p_s(\bz^k)-d_s(\blmbd^k;\bz^k)]\) to both sides leads to
\begin{equation}
\label{eq:b2-prime-dual-nesterov}
\begin{aligned}
    &p_s(\bz^k)-d_s(\blmbd^{k+1};\bz^k)
    +\frac{1}{2s}\|\bw^{k+1}-\blmbd_s^\star(\bz^k)\|^2\\
    &\quad-\Bigl(1-\frac{1}{\sqrt{\kappa}}\Bigr)
    \Bigl[p_s(\bz^k)-d_s(\blmbd^k;\bz^k)
    +\frac{1}{2s}\|\bw^k-\blmbd_s^\star(\bz^k)\|^2\Bigr]\\
    &\le-\frac{1-\tau L_\lambda}{2\tau}\|\blmbd^{k+1}-\bar{\blmbd}^k\|^2
    +\bigl[\mathcal P_s(\bx^{k+1},\blmbd^{k+1};\bz^k)
    -d_s(\blmbd^{k+1};\bz^k)\bigr]\\
    &\quad-\Bigl(1-\frac{1}{\sqrt{\kappa}}\Bigr)
    \bigl[\mathcal P_s(\bx^{k+1},\blmbd^k;\bz^k)-d_s(\blmbd^k;\bz^k)\bigr]
    -\frac{1}{\sqrt{\kappa}}
    \bigl[\mathcal P_s(\bx^{k+1},\blmbd_s^\star(\bz^k);\bz^k)-p_s(\bz^k)\bigr].
\end{aligned}
\end{equation}

\emph{Step 2 (Primal-subproblem error and absorption).} It remains to control
the three brackets on the right-hand side of~\eqref{eq:b2-prime-dual-nesterov},
which requires \(\bx^{k+1}\) to approximate \(\bx_s(\bar{\blmbd}^k;\bz^k)\). Under~\eqref{eq:b2-prime-stepsizes},
the primal step of Algorithm~\ref{alg:b2-prime} is a gradient descent step
with stepsize \(\gamma=2/(\mu_\gamma+L_x)\) on the \(\mu_\gamma\)-strongly
convex and \(L_x\)-smooth function
\(\mathcal P_s(\cdot,\bar{\blmbd}^k;\bz^k)\), whose minimizer is
\(\bx_s(\bar{\blmbd}^k;\bz^k)\), so the contraction of gradient
descent~\cite{nesterov2018lectures} and \(L_x\)-smoothness give
\begin{align*}
    \|\bx^{k+1}-\bx_s(\bar{\blmbd}^k;\bz^k)\|
    &\le\frac{L_x-\mu_\gamma}{L_x+\mu_\gamma}
    \|\bx^k-\bx_s(\bar{\blmbd}^k;\bz^k)\|
    =\gamma L_{\calD,\bx}\|\bx^k-\bx_s(\bar{\blmbd}^k;\bz^k)\|,\\
    \|\bx^{k+1}-\bx^k\|
    &\le\gamma L_x\|\bx^k-\bx_s(\bar{\blmbd}^k;\bz^k)\|.
\end{align*}
Since \(\bx_s(\bar{\blmbd}^k;\bz^k)\) shifts with \(\bar{\blmbd}^k\), the triangle
inequality and Lemma~\ref{lem:b2-prime-fixed-center-values}(i) yield
\[
    \|\bx^k-\bx_s(\bar{\blmbd}^k;\bz^k)\|
    \le\|\bx^k-\bx_s(\bar{\blmbd}^{k-1};\bz^k)\|
    +\frac{G_{\calA}}{\mu_\gamma}\|\blmbd^k-\bar{\blmbd}^{k-1}\|
    +\|\bx_s(\blmbd^k;\bz^k)-\bx_s(\bar{\blmbd}^k;\bz^k)\|.
\]
Let
\[
    S_k:=\|\bx^k-\bx_s(\bar{\blmbd}^{k-1};\bz^k)\|^2
    +\frac{G_{\calA}^2}{\mu_\gamma^2}\|\blmbd^k-\bar{\blmbd}^{k-1}\|^2
    +\|\bx_s(\bar{\blmbd}^k;\bz^k)-\bx_s(\blmbd^k;\bz^k)\|^2 .
\]
Combining the last three inequalities with
\((a+b+c)^2\le3(a^2+b^2+c^2)\) and~\eqref{eq:b2-prime-stepsizes}, we obtain
\begin{equation}
\label{eq:b2-prime-tracking}
    \|\bx^{k+1}-\bx_s(\bar{\blmbd}^k;\bz^k)\|^2
    \le3(\gamma L_{\calD,\bx})^2S_k\le\frac{3}{100}S_k,
    \qquad
    \|\bx^{k+1}-\bx^k\|^2
    \le3(\gamma L_x)^2S_k\le\frac{363}{100}S_k .
\end{equation}
With this primal-error estimate, Lemma~\ref{lem:b2-prime-fixed-center-values}(i)--(ii)
and the optimality of \(\bx_s^\star(\bz^k)\) for
\(\mathcal P_s(\cdot,\blmbd_s^\star(\bz^k);\bz^k)\) bound the three brackets as
\begin{align*}
    \mathcal P_s(\bx^{k+1},\blmbd^{k+1};\bz^k)-d_s(\blmbd^{k+1};\bz^k)
    &\le L_x\|\bx^{k+1}-\bx_s(\bar{\blmbd}^k;\bz^k)\|^2
    +\frac{L_xG_{\calA}^2}{\mu_\gamma^2}\|\blmbd^{k+1}-\bar{\blmbd}^k\|^2,\\
    \mathcal P_s(\bx^{k+1},\blmbd^k;\bz^k)-d_s(\blmbd^k;\bz^k)
    &\ge\frac{\mu_\gamma}{4}
    \|\bx_s(\bar{\blmbd}^k;\bz^k)-\bx_s(\blmbd^k;\bz^k)\|^2
    -\frac{\mu_\gamma}{2}\|\bx^{k+1}-\bx_s(\bar{\blmbd}^k;\bz^k)\|^2,\\
    \mathcal P_s(\bx^{k+1},\blmbd_s^\star(\bz^k);\bz^k)-p_s(\bz^k)
    &\ge\frac{\mu_\gamma}{2}\|\bx^{k+1}-\bx_s^\star(\bz^k)\|^2 .
\end{align*}
Substituting these bounds into~\eqref{eq:b2-prime-dual-nesterov},
adding the last two terms of
\(\mathcal E_{k+1}(\bz^k)-(1-\kappa^{-1/2})\mathcal E_k\), and noting
\(\frac12\le1-\kappa^{-1/2}\le1\) since \(\kappa=s/\tau\ge4\), we find
\[
\begin{aligned}
    \mathcal E_{k+1}(\bz^k)-\Bigl(1-\frac{1}{\sqrt{\kappa}}\Bigr)\mathcal E_k
    \le{}&-\left[\frac{1-\tau L_\lambda}{2\tau}
    -\frac{G_{\calA}^2}{\mu_\gamma^2}
    \Bigl(L_x+\frac{\mu_\gamma}{4}\Bigr)\right]
    \|\blmbd^{k+1}-\bar{\blmbd}^k\|^2\\
    &+\Bigl(L_x+\frac{3\mu_\gamma}{4}\Bigr)
    \|\bx^{k+1}-\bx_s(\bar{\blmbd}^k;\bz^k)\|^2
    -\frac{\mu_\gamma}{8}S_k
    -\frac{\mu_\gamma}{2\sqrt{\kappa}}\|\bx^{k+1}-\bx_s^\star(\bz^k)\|^2 .
\end{aligned}
\]
Finally, we bound the two coefficients on the right-hand side by
\begin{align*}
    \frac{\mu_\gamma}{8}S_k
    -\Bigl(L_x+\frac{3\mu_\gamma}{4}\Bigr)
    \|\bx^{k+1}-\bx_s(\bar{\blmbd}^k;\bz^k)\|^2
    &\ge\frac{57}{1000\gamma}S_k
    \ge\frac{1}{64\gamma}\|\bx^{k+1}-\bx^k\|^2,\\
    \frac{1-\tau L_\lambda}{2\tau}
    -\frac{G_{\calA}^2}{\mu_\gamma^2}\Bigl(L_x+\frac{\mu_\gamma}{4}\Bigr)
    &\ge\frac{1}{8\tau},
\end{align*}
where the first chain uses~\eqref{eq:b2-prime-tracking} together with
\(\mu_\gamma\ge\frac{9}{10\gamma}\) and \(L_x\le\frac{11}{10\gamma}\)
under~\eqref{eq:b2-prime-stepsizes}, and the second uses
\(\tau L_\lambda\le\frac12\) and \(14\gamma\tau G_{\calA}^2\le1\)
from~\eqref{eq:b2-prime-stepsizes}. Substituting both bounds into the preceding
display yields the claimed inequality and completes the proof.
\end{proof}

We next bound the proximal-descent term in~\eqref{eq:b2-prime-decomposition}.

\begin{lemma}
\label{lem:b2-prime-z-step-p}
Suppose that Assumptions~\ref{ass:smoothness} and~\ref{ass:h} hold, and let
\(0\le\alpha\le1\). Then, for every \(k\ge0\),
\[
    p_s(\bz^{k+1})-p_s(\bz^k)
    \le-\frac{\alpha(1-\alpha)}{2\gamma}\|\bx^{k+1}-\bz^k\|^2
    +\frac{\alpha}{2\gamma}\|\bx^{k+1}-\bx_s^\star(\bz^k)\|^2 .
\]
\end{lemma}

\begin{proof}
The center step of Algorithm~\ref{alg:b2-prime} gives
\(\bz^{k+1}-\bz^k=\alpha(\bx^{k+1}-\bz^k)\), so applying the first
inequality of Lemma~\ref{lem:b2-prime-center-dependence} with \(\bz=\bz^k\)
and \(\bz'=\bz^{k+1}\) yields
\begin{align*}
    p_s(\bz^{k+1})-p_s(\bz^k)
    &\le\frac{\alpha}{\gamma}\bigl\langle\bz^k-\bx_s^\star(\bz^k),
    \bx^{k+1}-\bz^k\bigr\rangle
    +\frac{\alpha^2}{2\gamma}\|\bx^{k+1}-\bz^k\|^2\\
    &=-\frac{\alpha(2-\alpha)}{2\gamma}\|\bx^{k+1}-\bz^k\|^2
    +\frac{\alpha}{\gamma}\bigl\langle\bx^{k+1}-\bx_s^\star(\bz^k),
    \bx^{k+1}-\bz^k\bigr\rangle\\
    &\le-\frac{\alpha(1-\alpha)}{2\gamma}\|\bx^{k+1}-\bz^k\|^2
    +\frac{\alpha}{2\gamma}\|\bx^{k+1}-\bx_s^\star(\bz^k)\|^2,
\end{align*}
where the equality splits
\(\bz^k-\bx_s^\star(\bz^k)
=(\bx^{k+1}-\bx_s^\star(\bz^k))-(\bx^{k+1}-\bz^k)\), and the last step
follows from Young's inequality, which completes the proof.
\end{proof}

We next carry out Step~3. The following lemma bounds the change of
\(\mathcal E_{k+1}\) caused by shifting the proximal center.

\begin{lemma}
\label{lem:b2-prime-center-error}
Suppose that Assumptions~\ref{ass:smoothness} and~\ref{ass:h} hold. Then,
for every \(k\ge0\),
\[
    \mathcal E_{k+1}(\bz^{k+1})-\mathcal E_{k+1}(\bz^k)
    \le\frac{\alpha}{\gamma}\sqrt{\frac{7}{2\mu_\gamma}}
    \sqrt{\mathcal E_{k+1}(\bz^k)}\,\|\bx^{k+1}-\bz^k\|
    +\frac{11\alpha^2}{8\gamma^2\mu_\gamma}\|\bx^{k+1}-\bz^k\|^2 .
\]
\end{lemma}

\begin{proof}
The fourth term of \(\mathcal E_{k+1}(\bz)\) is independent of \(\bz\), so
we bound the change of the first three terms from \(\bz^k\) to \(\bz^{k+1}\).
For the primal-dual gap, Lemmas~\ref{lem:b2-prime-fixed-center-values}
and~\ref{lem:b2-prime-center-dependence} give
\begin{align*}
    &\bigl[p_s(\bz^{k+1})-d_s(\blmbd^{k+1};\bz^{k+1})\bigr]
    -\bigl[p_s(\bz^k)-d_s(\blmbd^{k+1};\bz^k)\bigr]\\
    &\qquad\le\frac{1}{\gamma}\bigl\langle
    \bx_s(\blmbd^{k+1};\bz^{k+1})-\bx_s^\star(\bz^k),\bz^{k+1}-\bz^k\bigr\rangle\\
    &\qquad\le\frac{1}{\gamma}\|\bx_s(\blmbd^{k+1};\bz^k)-\bx_s^\star(\bz^k)\|\,
    \|\bz^{k+1}-\bz^k\|
    +\frac{1}{\gamma^2\mu_\gamma}\|\bz^{k+1}-\bz^k\|^2\\
    &\qquad\le\frac{\sqrt2}{\gamma\sqrt{\mu_\gamma}}
    \sqrt{p_s(\bz^k)-d_s(\blmbd^{k+1};\bz^k)}\,\|\bz^{k+1}-\bz^k\|
    +\frac{1}{\gamma^2\mu_\gamma}\|\bz^{k+1}-\bz^k\|^2,
\end{align*}
where the first inequality subtracts the second bound of
Lemma~\ref{lem:b2-prime-center-dependence} from the first, the second uses the
Lipschitz continuity of \(\bx_s(\blmbd^{k+1};\cdot)\), and the last follows
from Lemma~\ref{lem:b2-prime-fixed-center-values}(iii). The same argument, now with the Lipschitz bounds on \(\blmbd_s^\star\) and
\(\bx_s(\bar{\blmbd}^k;\cdot)\), bounds the remaining two terms as
\begin{align*}
    &\frac{1}{2s}\|\bw^{k+1}-\blmbd_s^\star(\bz^{k+1})\|^2
    -\frac{1}{2s}\|\bw^{k+1}-\blmbd_s^\star(\bz^k)\|^2\\
    &\qquad\le\frac{1}{\gamma\sqrt{2\mu_\gamma}}
    \sqrt{\frac{1}{2s}\|\bw^{k+1}-\blmbd_s^\star(\bz^k)\|^2}\,
    \|\bz^{k+1}-\bz^k\|
    +\frac{1}{8\gamma^2\mu_\gamma}\|\bz^{k+1}-\bz^k\|^2,\\
    &\frac{\mu_\gamma}{4}\|\bx^{k+1}-\bx_s(\bar{\blmbd}^k;\bz^{k+1})\|^2
    -\frac{\mu_\gamma}{4}\|\bx^{k+1}-\bx_s(\bar{\blmbd}^k;\bz^k)\|^2\\
    &\qquad\le\frac{1}{\gamma\sqrt{\mu_\gamma}}
    \sqrt{\frac{\mu_\gamma}{4}\|\bx^{k+1}-\bx_s(\bar{\blmbd}^k;\bz^k)\|^2}\,
    \|\bz^{k+1}-\bz^k\|
    +\frac{1}{4\gamma^2\mu_\gamma}\|\bz^{k+1}-\bz^k\|^2 .
\end{align*}
Since the first three terms of \(\mathcal E_{k+1}(\bz^k)\) are nonnegative
and sum to at most \(\mathcal E_{k+1}(\bz^k)\), adding the three estimates and
applying the Cauchy--Schwarz inequality to the coefficients
\((\sqrt2,1/\sqrt2,1)\) of the linear terms gives
\[
    \mathcal E_{k+1}(\bz^{k+1})-\mathcal E_{k+1}(\bz^k)
    \le\frac{1}{\gamma}\sqrt{\frac{7}{2\mu_\gamma}}
    \sqrt{\mathcal E_{k+1}(\bz^k)}\,\|\bz^{k+1}-\bz^k\|
    +\frac{11}{8\gamma^2\mu_\gamma}\|\bz^{k+1}-\bz^k\|^2 .
\]
Substituting \(\|\bz^{k+1}-\bz^k\|=\alpha\|\bx^{k+1}-\bz^k\|\), which follows
from the center step~\eqref{eq:fast-lralm-center-step}, completes the proof.
\end{proof}

Combining Lemmas~\ref{lem:b2-prime-fixed-center-contraction}
and~\ref{lem:b2-prime-z-step-p} from Step~2 with
Lemma~\ref{lem:b2-prime-center-error} from Step~3 yields the sufficient descent
and an \(s\)-independent initial gap for \(\{\mathcal W_k\}_{k\ge0}\).

\begin{lemma}
\label{lem:b2-prime-potential-descent}
Suppose that Assumptions~\ref{ass:smoothness} and~\ref{ass:h} and
\eqref{eq:b2-prime-stepsizes} hold. Then, for every \(k\ge0\),
\[
    \mathcal W_{k+1}-\mathcal W_k
    \le -\frac{1}{128\gamma}\|\bx^{k+1}-\bx^k\|^2
        -\frac{1}{16\tau}\|\blmbd^{k+1}-\bar{\blmbd}^k\|^2
        -\frac{\alpha}{6\gamma}\|\bx^{k+1}-\bz^k\|^2.
\]
Moreover, there exists a constant \(\tilde{\Delta}>0\), depending only on
\(\mathbf C\), \(\Delta\), \(\rho\), \(\gamma\), and \(\tau\), such that for
every \(K\ge0\),
\begin{equation}
    \mathcal W_0-\mathcal W_K\le\tilde{\Delta}.
    \label{eq:b2-prime-initial-gap}
\end{equation}
\end{lemma}

\begin{proof}
Lemma~\ref{lem:b2-prime-center-error} combined with Young's inequality
gives
\[
    \mathcal E_{k+1}
    \le\Bigl(1+\frac{1}{2\sqrt{\kappa}}\Bigr)\mathcal E_{k+1}(\bz^k)
    +\frac{14\sqrt{\kappa}+11}{8\gamma^2\mu_\gamma}
    \alpha^2\|\bx^{k+1}-\bz^k\|^2 .
\]
We next insert Lemma~\ref{lem:b2-prime-fixed-center-contraction} into this
bound and apply
\((1+\frac{1}{2\sqrt{\kappa}})(1-\frac{1}{\sqrt{\kappa}})
\le1-\frac{1}{2\sqrt{\kappa}}\). Adding half of the resulting estimate on
\(\mathcal E_{k+1}-\mathcal E_k\) to Lemma~\ref{lem:b2-prime-z-step-p}
via~\eqref{eq:b2-prime-decomposition} and dropping the nonpositive term
\(-\frac{1}{4\sqrt{\kappa}}\mathcal E_k\) yields
\[
\begin{aligned}
    \mathcal W_{k+1}-\mathcal W_k
    \le{}&-\frac{1}{128\gamma}\|\bx^{k+1}-\bx^k\|^2
    -\frac{1}{16\tau}\|\blmbd^{k+1}-\bar{\blmbd}^k\|^2
    +\Bigl(\frac{\alpha}{2\gamma}-\frac{\mu_\gamma}{4\sqrt{\kappa}}\Bigr)
    \|\bx^{k+1}-\bx_s^\star(\bz^k)\|^2\\
    &-\frac{\alpha}{\gamma}
    \Bigl(\frac{1-\alpha}{2}
    -\frac{14\sqrt{\kappa}+11}{16\gamma\mu_\gamma}\alpha\Bigr)
    \|\bx^{k+1}-\bz^k\|^2 .
\end{aligned}
\]
Under~\eqref{eq:b2-prime-stepsizes}, the coefficient of
\(\|\bx^{k+1}-\bx_s^\star(\bz^k)\|^2\) is nonpositive and that of
\(\|\bx^{k+1}-\bz^k\|^2\) is at most \(-\frac{\alpha}{6\gamma}\), which proves
the descent inequality.

To prove~\eqref{eq:b2-prime-initial-gap}, we bound \(\mathcal W_K\) from below
and \(\mathcal W_0\) from above. By~\eqref{eq:reduced-al}, for
any \(\sigma>0\) and \(\bx\in\RR^n\), we have
\(\calD_\sigma(\bx,\zero)=f(\bx)+\min_{\by\in\RR^m}\{h(\by)
+\frac{\sigma}{2}\|\calA(\bx)-\by\|^2\}\). Choosing \(\by=\calA(\bx)\) in the minimum and using
\(h(\by)\ge h(\calA(\bx))-L_h\|\calA(\bx)-\by\|\) give
\begin{equation}
    F(\bx)-\frac{L_h^2}{2\sigma}\le\calD_\sigma(\bx,\zero)\le F(\bx).
    \label{eq:b2-prime-D-bounds}
\end{equation}
By~\eqref{eq:minimax-equality} and~\eqref{eq:effective-penalty} with
\(\bmu=\zero\), we have
\(p_s(\bz)=\min_{\bx\in\RR^n}\{\calD_{\rho+s}(\bx,\zero)
+\frac{1}{2\gamma}\|\bx-\bz\|^2\}\), and~\eqref{eq:b2-prime-D-bounds} with
\(\sigma=\rho+s\) yields
\(\inf F-\frac{L_h^2}{2(\rho+s)}\le p_s(\bz)\le F(\bz)\) for all
\(\bz\in\RR^n\). As \(p_s(\bz^K)\ge d_s(\blmbd^K;\bz^K)\) gives \(\mathcal E_K\ge0\), it
follows that
\begin{equation}
    \mathcal W_K=p_s(\bz^K)+\frac12\mathcal E_K
    \ge\inf_{\bx\in\RR^n}F(\bx)-\frac{L_h^2}{2(\rho+s)} .
    \label{eq:b2-prime-WK-lower}
\end{equation}
For \(\mathcal W_0\), the initialization conditions \(\bx^0=\zero\),
\(\|\blmbd^0\|\le L_h\), \(\|\bz^0\|\le1\), and \(\blmbd^{-1}=\blmbd^0\) give
\(\bar{\blmbd}^{-1}=\bw^0=\blmbd^0\), so
\begin{equation}
    \mathcal W_0=p_s(\bz^0)+\frac12\Bigl[
    p_s(\bz^0)-d_s(\blmbd^0;\bz^0)
    +\frac{1}{2s}\|\blmbd^0-\blmbd_s^\star(\bz^0)\|^2
    +\frac{\mu_\gamma}{4}\|\bx_s(\blmbd^0;\bz^0)\|^2\Bigr].
    \label{eq:b2-prime-W0}
\end{equation}
Set \(G_0:=G_f+G_{\calA}L_h+1/\gamma\). By Lemma~\ref{lem:reduced-al}(a),
Assumptions~\ref{ass:smoothness} and~\ref{ass:h}, and \(\|\bz^0\|\le1\), we
have \(\|\bq(\zero,\blmbd^0)\|\le L_h\) and
\[
    \|\nabla_{\bx}\mathcal P_s(\zero,\blmbd^0;\bz^0)\|
    =\Bigl\|\nabla_{\bx}\calD_\rho(\zero,\blmbd^0)-\frac{1}{\gamma}\bz^0\Bigr\|
    \le G_0.
\]
Next, for every \(\by\in\RR^m\), the \(L_h\)-Lipschitz continuity of \(h\) and
the Cauchy--Schwarz inequality with \(\|\blmbd^0\|\le L_h\) yield
\begin{align*}
    h(\by)+\langle\blmbd^0,\calA(\zero)-\by\rangle
    +\frac{\rho}{2}\|\calA(\zero)-\by\|^2
    &\ge h(\calA(\zero))-2L_h\|\calA(\zero)-\by\|
    +\frac{\rho}{2}\|\calA(\zero)-\by\|^2\\
    &\ge h(\calA(\zero))-\frac{2L_h^2}{\rho},
\end{align*}
where the second inequality completes the square in \(\|\calA(\zero)-\by\|\).
Minimizing over \(\by\in\RR^m\) in~\eqref{eq:reduced-al} gives
\(\calD_\rho(\zero,\blmbd^0)\ge F(\zero)-\frac{2L_h^2}{\rho}\). Together with
\(\|\blmbd^0\|\le L_h\) in~\eqref{eq:potential}, we obtain
\[
    \mathcal P_s(\zero,\blmbd^0;\bz^0)
    =\calD_\rho(\zero,\blmbd^0)+\frac{1}{2\gamma}\|\bz^0\|^2
    -\frac{1}{2s}\|\blmbd^0\|^2
    \ge F(\zero)-\frac{2L_h^2}{\rho}-\frac{L_h^2}{2s}.
\]
By the \(\mu_\gamma\)-strong convexity of
\(\mathcal P_s(\cdot,\blmbd^0;\bz^0)\), its minimizer \(\bx_s(\blmbd^0;\bz^0)\)
and minimum value \(d_s(\blmbd^0;\bz^0)\) satisfy
\begin{align*}
    \|\bx_s(\blmbd^0;\bz^0)\|
    &\le\frac{1}{\mu_\gamma}\|\nabla_{\bx}\mathcal P_s(\zero,\blmbd^0;\bz^0)\|
    \le\frac{G_0}{\mu_\gamma},\\
    d_s(\blmbd^0;\bz^0)
    &\ge\mathcal P_s(\zero,\blmbd^0;\bz^0)
    -\frac{1}{2\mu_\gamma}\|\nabla_{\bx}\mathcal P_s(\zero,\blmbd^0;\bz^0)\|^2
    \ge F(\zero)-\frac{2L_h^2}{\rho}-\frac{L_h^2}{2s}-\frac{G_0^2}{2\mu_\gamma}.
\end{align*}
Moreover, applying~\eqref{eq:effective-penalty-maximizer} at
\(\bx=\bx_s^\star(\bz^0)\) and \(\bmu=\zero\) gives
\(\blmbd_s^\star(\bz^0)
=\frac{s}{\rho+s}\bq(\bx_s^\star(\bz^0),\blmbd_s^\star(\bz^0))\), so
\(\|\blmbd_s^\star(\bz^0)\|\le L_h\) by Lemma~\ref{lem:reduced-al}(a), and
using \(\|\blmbd^0\|\le L_h\) yields
\[
    \frac{1}{2s}\|\blmbd^0-\blmbd_s^\star(\bz^0)\|^2\le\frac{2L_h^2}{s}.
\]
Evaluating \(p_s(\bz^0)\) at \(\bx=\zero\) and applying~\eqref{eq:b2-prime-D-bounds}
with \(\sigma=\rho+s\) and \(\|\bz^0\|\le1\) give
\[
    p_s(\bz^0)\le\calD_{\rho+s}(\zero,\zero)+\frac{1}{2\gamma}\|\bz^0\|^2
    \le F(\zero)+\frac{1}{2\gamma}.
\]
Substituting these four estimates into~\eqref{eq:b2-prime-W0} and subtracting~\eqref{eq:b2-prime-WK-lower}, we obtain
\begin{align*}
    \mathcal W_0-\mathcal W_K
    &\le\bigl(F(\zero)-\inf_{\bx\in\RR^n}F(\bx)\bigr)
    +\frac{3}{4\gamma}+\frac{L_h^2}{\rho}+\frac{5L_h^2}{4s}
    +\frac{L_h^2}{2(\rho+s)}+\frac{3G_0^2}{8\mu_\gamma}\\
    &\le\Delta+\frac1\gamma+\frac{2L_h^2}{\rho}+\frac{L_h^2}{2\tau}
    +\frac{G_0^2}{2\mu_\gamma}=:\tilde{\Delta},
\end{align*}
where the second inequality uses \(F(\zero)-\inf F\le\Delta\),
\(\rho+s\ge\rho\), and \(s\ge4\tau\) from~\eqref{eq:b2-prime-stepsizes}. Since
\(\tilde{\Delta}\) depends only on \(\mathbf C\), \(\Delta\), \(\rho\),
\(\gamma\), and \(\tau\), this establishes~\eqref{eq:b2-prime-initial-gap} and
completes the proof.
\end{proof}

Combining Lemma~\ref{lem:b2-prime-potential-descent} with the residual and
multiplier bounds in Lemmas~\ref{lem:b2-prime-residual-decomposition}
and~\ref{lem:b2-prime-multiplier-bounds} yields the iteration complexity of
Fast LRALM.

\begin{theorem}
\label{thm:b2-prime-complexity}
Suppose that Assumptions~\ref{ass:smoothness} and~\ref{ass:h} hold, and
let \(\rho,\gamma,\tau,c_\alpha,c_s>0\) be independent of \(I\), \(\bu^0\), and
\(\epsilon\) with
\(c_s\ge12\sqrt{\rho/\tau}\,L_h\). Given \(0<\epsilon\le c_s/(4\tau)\),
let \(\{(\bx^k,\blmbd^k,\bz^k)\}_{k\ge0}\) be generated by
Algorithm~\ref{alg:b2-prime} under~\eqref{eq:b2-prime-stepsizes} with
\(s=c_s/\epsilon\). Then,
for every \(K\ge1\),
\begin{equation}
    \min_{0\le k<K}\mathcal R_{\mathrm{KKT}}(\bx^{k+1},\by^{k+1},\bq^k)
    \le\frac{C_R\sqrt{\sqrt{\kappa}\,\tilde{\Delta}}}{\sqrt{K}}
    +\frac{6\sqrt{\rho/\tau}\,L_h}{s},
    \label{eq:b2-prime-rate}
\end{equation}
where \(\by^{k+1}\) and \(\bq^k\) are defined in~\eqref{eq:b2-prime-outputs},
\(C_R:=\max\bigl\{8\sqrt{\gamma}\,L_{\calD,\bx}+\sqrt{6/(c_\alpha\gamma)},\,
\sqrt{8/\tau}\bigr\}\), and \(\tilde{\Delta}\) is the constant in
Lemma~\ref{lem:b2-prime-potential-descent}. Consequently, the iteration complexity of Fast LRALM to obtain an \(\epsilon\)-KKT point is of order
\(\mathcal O(\epsilon^{-5/2})\).
\end{theorem}
\begin{proof}
Summing the descent inequality in
Lemma~\ref{lem:b2-prime-potential-descent} from \(k=0\) to \(K-1\) and
dividing by \(K\), together with~\eqref{eq:b2-prime-initial-gap} and
\(\alpha=\frac{c_\alpha}{\sqrt{\kappa}}\), yields
\[
    \frac1K\sum_{k=0}^{K-1}\Bigl(
    \frac{1}{128\gamma}\|\bx^{k+1}-\bx^k\|^2
    +\frac{1}{16\tau}\|\blmbd^{k+1}-\bar{\blmbd}^k\|^2
    +\frac{c_\alpha}{6\gamma\sqrt{\kappa}}\|\bx^{k+1}-\bz^k\|^2\Bigr)
    \le\frac{\tilde{\Delta}}{K}.
\]
Hence there exists \(0\le k<K\) such that
\[
    \|\bx^{k+1}-\bx^k\|\le8\sqrt{\frac{2\gamma\tilde{\Delta}}{K}},
    \qquad
    \|\blmbd^{k+1}-\bar{\blmbd}^k\|\le4\sqrt{\frac{\tau\tilde{\Delta}}{K}},
    \qquad
    \|\bx^{k+1}-\bz^k\|\le\sqrt{\frac{6\gamma\sqrt{\kappa}\,\tilde{\Delta}}{c_\alpha K}}.
\]
Because~\eqref{eq:b2-prime-stepsizes} ensures \(\tau L_\lambda\le\frac12\),
Lemma~\ref{lem:b2-prime-residual-decomposition} with \(\eta=\gamma\) and
Lemma~\ref{lem:b2-prime-multiplier-bounds} give
\[
    \mathcal R_{\mathrm{KKT}}(\bx^{k+1},\by^{k+1},\bq^k)
    \le\max\Bigl\{L_{\calD,\bx}\|\bx^{k+1}-\bx^k\|
    +\frac1\gamma\|\bx^{k+1}-\bz^k\|,\
    \frac1\tau\|\blmbd^{k+1}-\bar{\blmbd}^k\|\Bigr\}
    +\frac{6\sqrt{\rho/\tau}\,L_h}{s}.
\]
Substituting the three bounds above and using \(\kappa=s/\tau\ge4\)
yield~\eqref{eq:b2-prime-rate}.

Finally, with \(s=c_s/\epsilon\) and \(c_s\ge12\sqrt{\rho/\tau}\,L_h\), the
second term of~\eqref{eq:b2-prime-rate} is at most \(\epsilon/2\), and the
first term is at most \(\epsilon/2\) whenever
\(K\ge4C_R^2\tilde{\Delta}\sqrt{c_s/\tau}\,\epsilon^{-5/2}\), which
gives \(\min_{0\le k<K}\mathcal R_{\mathrm{KKT}}(\bx^{k+1},\by^{k+1},\bq^k)\le\epsilon\).
Since \(C_R\) and \(\tilde{\Delta}\) are independent of \(\epsilon\), the
iteration complexity is \(\mathcal O(\epsilon^{-5/2})\), which completes the
proof.
\end{proof}

\begin{remark}
\begin{enumerate}[label=(\roman*)]
\item Recent minimax methods also achieve \(\mathcal O(\epsilon^{-5/2})\)
complexity on smooth NC-C problems, including
the multi-loop methods~\cite{kong2021accelerated,ke2026efficient,pan2026optimal} and the single-loop projected extragradient method MC-EGM of Zhang
et~al.~\cite{zhang2026matching}. Complementing
these results, Fast LRALM operates on the unconstrained smooth dual of
\(\calD_\rho\) and uses the uniform bound on the updated multiplier \(\bq\) to
achieve the same \(\mathcal O(\epsilon^{-5/2})\) rate in a single loop without
dual projections or extragradient steps, using one call to \(\mathcal O_{f,\calA}\)
per iteration.
\item In augmented Lagrangian methods, Nesterov acceleration has been studied
for convex constrained optimization~\cite{he2010acceleration,he2026augmented}.
Alongside these developments, Theorem~\ref{thm:b2-prime-complexity} extends
Nesterov multiplier extrapolation to the nonconvex composite
Problem~\eqref{prob:main} with a constant penalty \(\rho\),
establishing a single-loop \(\mathcal O(\epsilon^{-5/2})\) complexity guarantee.
Within the linearized ALM framework, this improves on the
\(\Theta(\epsilon^{-3})\) rate of \(\calC_\theta\) while keeping \(\rho\)
constant.
\end{enumerate}
\end{remark}

\section{Numerical Experiments}
\label{sec:experiments}
\setcounter{topnumber}{2}\setcounter{totalnumber}{3}\renewcommand{\topfraction}{0.9}\renewcommand{\textfraction}{0.07}\renewcommand{\floatpagefraction}{0.8}

We compare Fast LRALM in Algorithm~\ref{alg:b2-prime} on a nonconvex composite problem with three single-loop baselines, namely Smoothed GDA~\cite{zhang2020smoothed} and Perturbed Smoothed GDA~\cite{li2026smoothing} from \(\calC_\theta\) and the projected extragradient method MC-EGM~\cite{zhang2026matching}, together with the multi-loop inexact augmented Lagrangian method iALM (Sahin et~al.)~\cite{sahin2019inexact} with the inner solver of~\cite{luSingleLoopGradientDescent2022}.

\subsection{Experimental setup}
\label{subsec:exp-illcond}

We consider the nonconvex composite problem
\begin{equation}
  \min_{\bx\in\RR^n}\;
  \frac{1}{2}\sum_{l=1}^{n}\log\bigl(1+(x_l-d_l)^2\bigr)
  +\|\bB\bx-\bc\|_1,
  \label{eq:exp-illcond}
\end{equation}
which takes the form of Problem~\eqref{prob:main} with \(\calA(\bx)=\bB\bx-\bc\) and \(h=\|\cdot\|_1\). Specifically, \eqref{eq:exp-illcond} is a denoising variant of the Cauchy-loss model in Section~\ref{sec:introduction}. We set \(n=50\), \(m=20\), and \(\bB=\bU\operatorname{diag}(\bsigma)\bV^\top\in\RR^{m\times n}\), where \(\bU\in\RR^{m\times m}\) is an orthogonal matrix, \(\bV\in\RR^{n\times m}\) has orthonormal columns, and the entries of \(\bsigma\in\RR^m\) are logarithmically spaced from \(1\) down to \(\sigma_{\min}\in\{10^{-2},10^{-3},10^{-4}\}\). Since \(\bB\) has full row rank, \(\sigma_{\min}\) equals the smallest singular value of the Jacobian \(J_{\calA}(\bx)=\bB\); smaller values thus weaken the Jacobian regularity of \(\calA\) and increase ill-conditioning. The vectors \(\bc\) and \(\bd\) are constructed from a prescribed pair \((\bx^\star,\bq^\star)\) satisfying the KKT conditions of~\eqref{eq:exp-illcond} with \(\|\bq^\star\|_\infty<1\).

\begin{figure}[!t]
\centering
\includegraphics[width=\textwidth]{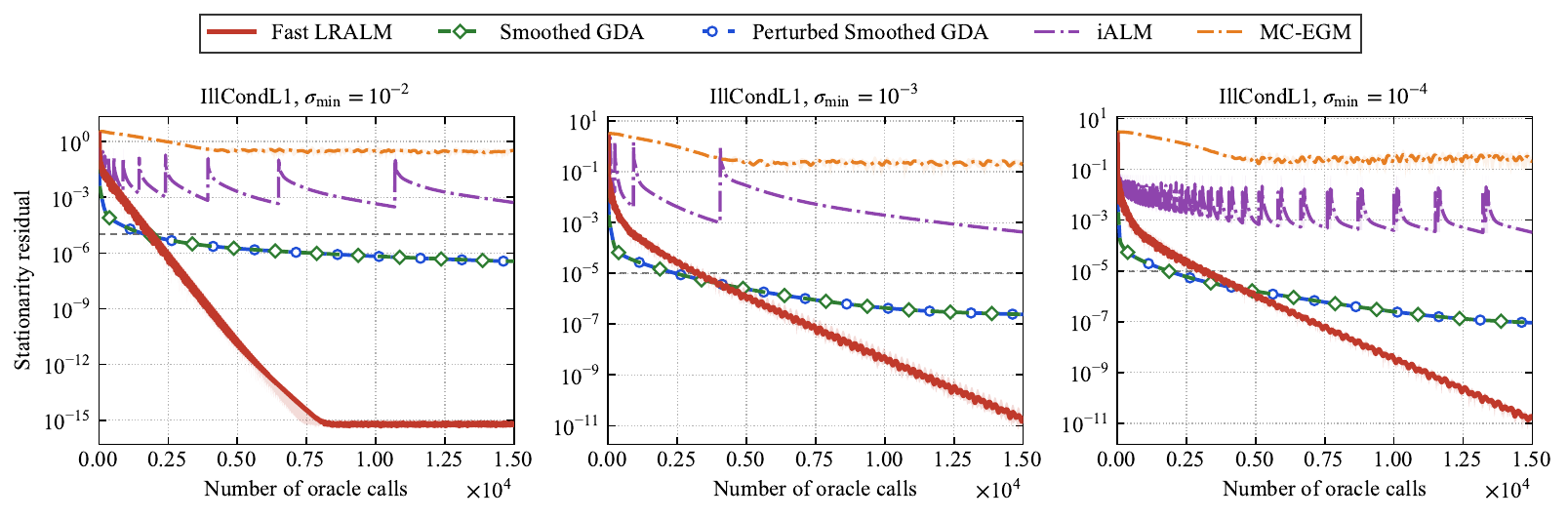}
\caption{Stationarity residual on Problem~\eqref{eq:exp-illcond}. The dashed line marks \(\epsilon=10^{-5}\).}
\label{fig:exp-stat}
\end{figure}

\begin{figure}[!t]
\centering
\includegraphics[width=\textwidth]{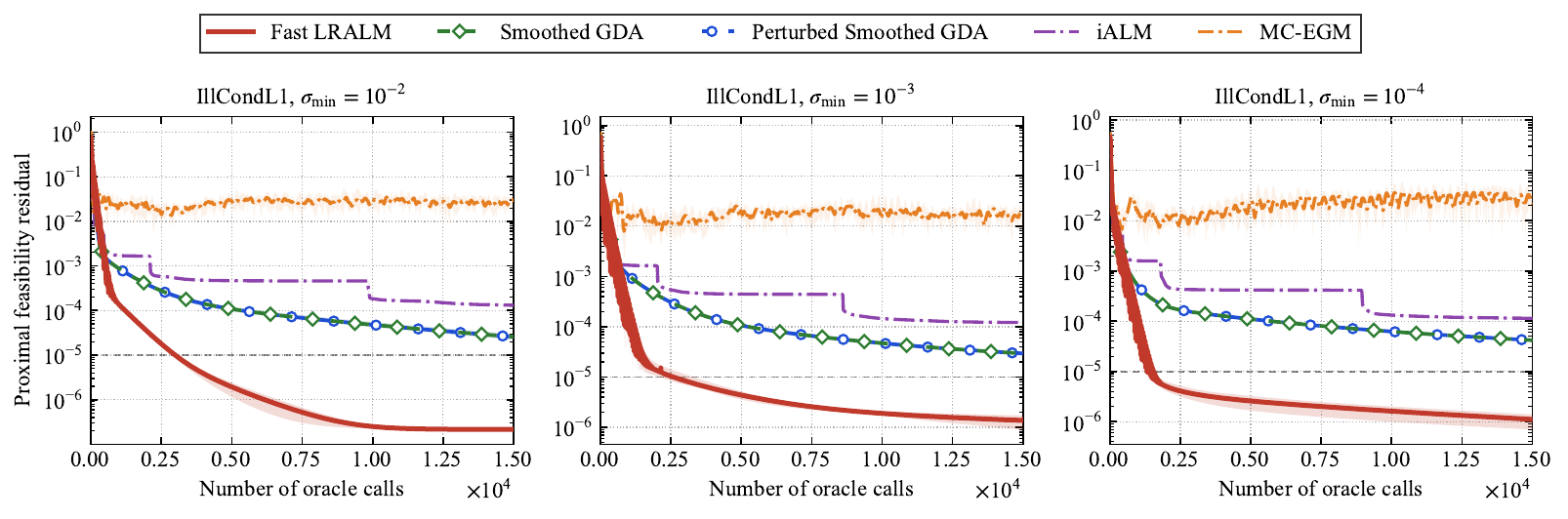}
\caption{Feasibility residual on Problem~\eqref{eq:exp-illcond}. The dashed line marks \(\epsilon=10^{-5}\).}
\label{fig:exp-feas}
\end{figure}

\begin{figure}[!t]
\centering
\includegraphics[width=\textwidth]{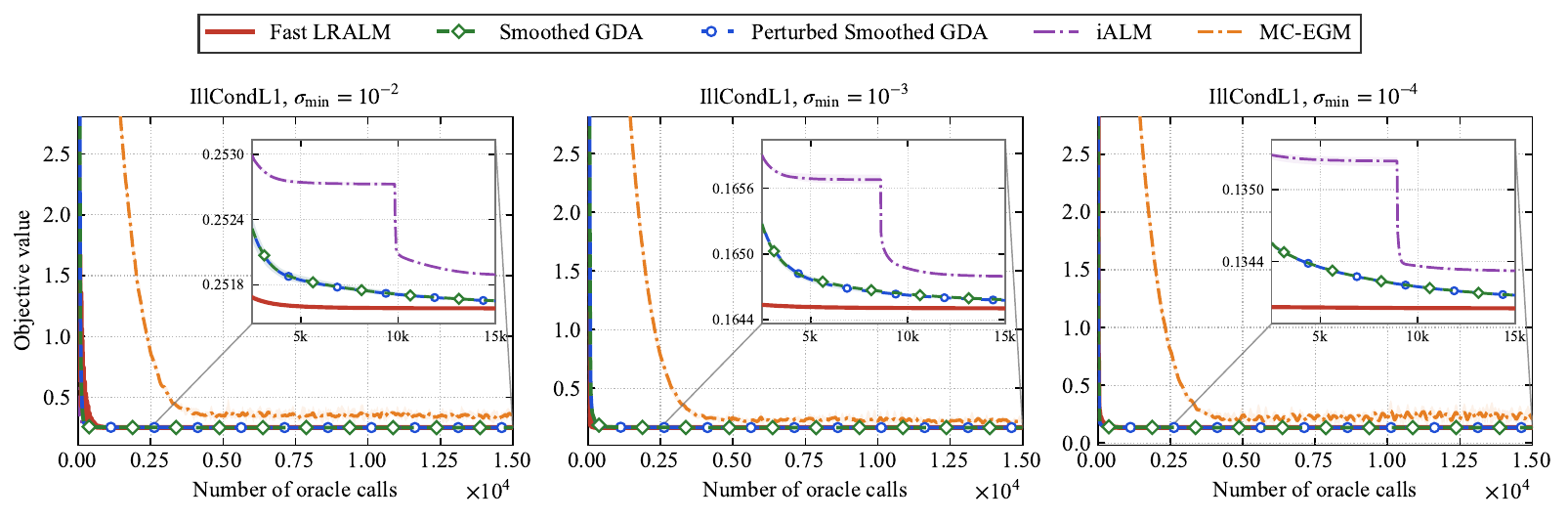}
\caption{Objective value on Problem~\eqref{eq:exp-illcond}. Insets zoom in on \(2{,}500\) to \(15{,}000\) oracle calls.}
\label{fig:exp-obj}
\end{figure}

All five methods run with \(\epsilon=10^{-5}\) under a shared budget of \(15{,}000\) calls to \(\mathcal{O}_{f,\calA}\), using one per iteration for Fast LRALM, Smoothed GDA, and Perturbed Smoothed GDA, two for MC-EGM, and all inner-loop calls for iALM. The first four methods solve~\eqref{eq:exp-illcond} via the reduced augmented Lagrangian \(\calD_\rho\); MC-EGM uses the equivalent smooth minimax formulation \(\min_{\bx\in\RR^n}\max_{\blmbd\in[-1,1]^m}\{f(\bx)+\langle\blmbd,\bB\bx-\bc\rangle\}\). Each curve in Figures~\ref{fig:exp-stat}--\ref{fig:exp-obj} plots the latest value averaged over five initial points, with shaded bands showing the min-to-max range.

\paragraph{Parameter selection.}
We tune hyperparameters for each method, instance, and metric by grid
search on the first initial point, selecting the configuration with
finite iterates and final objective at or below the initial value that
attains the smallest final metric value. Fast LRALM, Smoothed GDA, and
Perturbed Smoothed GDA share the grids \(\rho\in\{0.03,0.1,0.3,1\}\),
\(\gamma=c/L_{\calD,\bx}\) with \(\eta=\gamma\) and
\(c\in\{0.5,1,2,4\}\), \(\tau=t\rho\) with
\(t\in\{1/16,1/8,1/4,1/2\}\), and \(\alpha\in\{1/6,0.3,0.5,0.8\}\).
Smoothed GDA sets \(s=1\); Fast LRALM and Perturbed Smoothed GDA set
\(s=c_sL_h/\epsilon\) with \(c_s\in\{0.3,1,3,24\}\); Fast LRALM
additionally searches
\(\omega\in\{0.9,0.95,0.98,(\sqrt{\kappa}-1)/(\sqrt{\kappa}+1)\}\).
For the remaining baselines, iALM
follows~\cite{sahin2019inexact} with initial penalty \(1\), dual scale
\(10\), and four-point grids for its penalty growth factor and inner
stepsize, while MC-EGM~\cite{zhang2026matching} searches four-point
grids for its extragradient and center stepsizes.

\subsection{Numerical results}
\label{subsec:exp-results}

Figures~\ref{fig:exp-stat} and~\ref{fig:exp-feas} show the two components of \(\mathcal R_{\mathrm{KKT}}\) in~\eqref{eq:kkt-residual}. In Figure~\ref{fig:exp-stat}, Fast LRALM and both GDA methods drive the stationarity residual below \(\epsilon\) on all three instances, with Fast LRALM reaching several orders of magnitude lower. In Figure~\ref{fig:exp-feas}, only Fast LRALM brings the feasibility residual below \(\epsilon\) on every instance; both GDA methods plateau near \(\epsilon\), and iALM and MC-EGM remain above it. Smoothed GDA and Perturbed Smoothed GDA trace nearly identical curves throughout. In Figure~\ref{fig:exp-obj}, all methods except MC-EGM reach similar final objective values. Compared with the unaccelerated baselines, dual Nesterov extrapolation primarily accelerates feasibility reduction while preserving objective progress. This benefit persists across all three conditioning levels within a comparable oracle budget.

\begin{figure}[!t]
\centering
\includegraphics[width=0.8\textwidth]{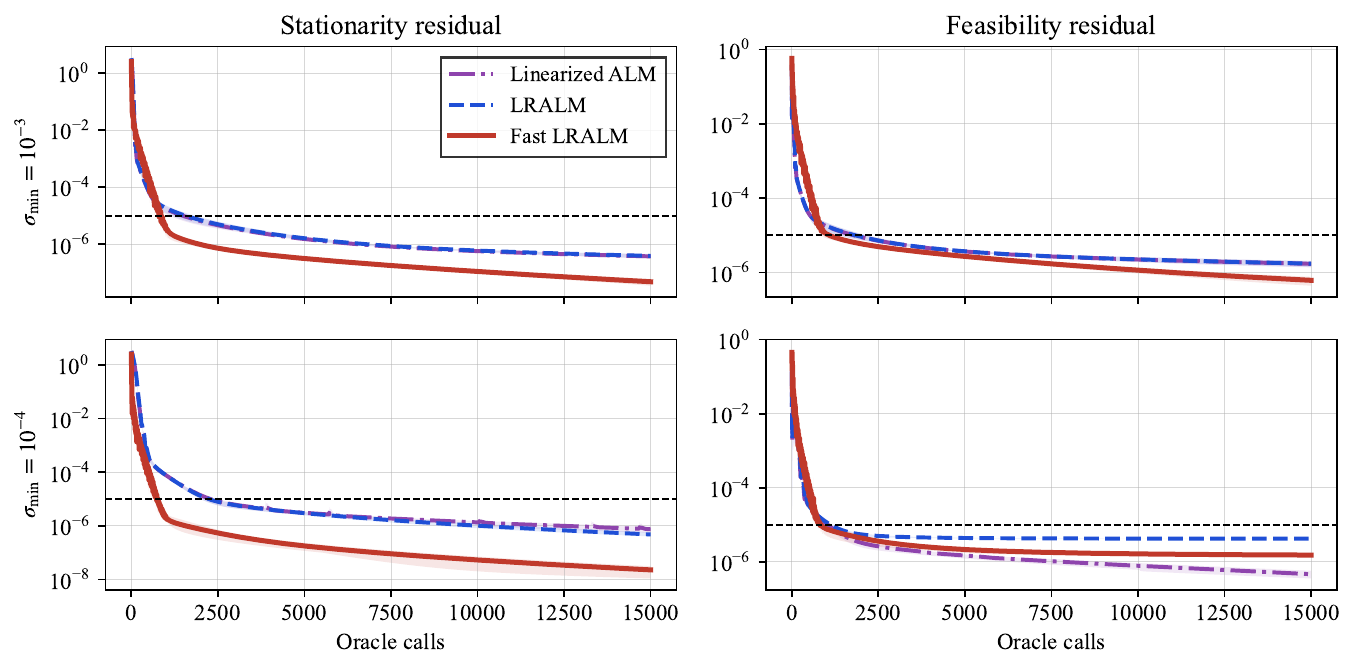}
\caption{Algorithm comparison on Problem~\eqref{eq:exp-illcond} under constant \(\rho\) with \(s=\Theta(\epsilon^{-1})\) for LRALM and Fast LRALM. The dashed line marks \(\epsilon=10^{-5}\).}
\label{fig:ablation-algos}
\end{figure}

\subsection{Ablation study}
\label{subsec:ablation}

To isolate the roles of \(\rho\), \(s\), and Nesterov dual extrapolation, we compare linearized ALM (\(\alpha=\beta=1\), \(s=\tau=\rho\)), LRALM (\(\beta=0\)), and Fast LRALM on~\eqref{eq:exp-illcond} for \(\sigma_{\min}\in\{10^{-3},10^{-4}\}\) under the same oracle budget. For linearized ALM we test constant \(\rho\) and \(\rho=\Theta(\epsilon^{-1})\); for LRALM and Fast LRALM we test three \((\rho,s)\) configurations, namely constant \(\rho\) with \(s=\Theta(\epsilon^{-1})\), constant \((\rho,s)\), and \(\rho,s=\Theta(\epsilon^{-1})\). Parameters are tuned by grid search.

\begin{figure}[!t]
\centering
\includegraphics[width=0.8\textwidth]{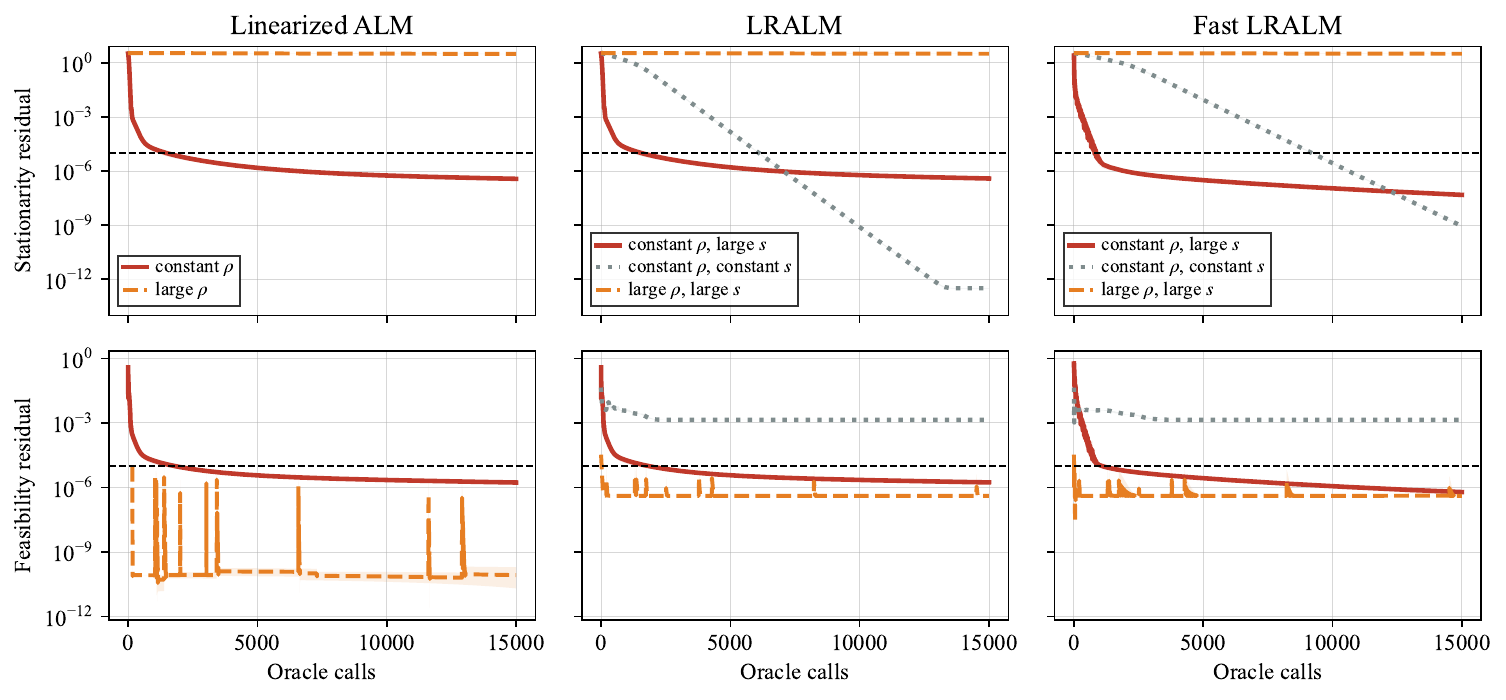}
\caption{Parameter ablation within each algorithm on Problem~\eqref{eq:exp-illcond} with \(\sigma_{\min}=10^{-3}\). The dashed line marks \(\epsilon=10^{-5}\).}
\label{fig:ablation-params-1e3}
\end{figure}

\begin{figure}[!t]
\centering
\includegraphics[width=0.8\textwidth]{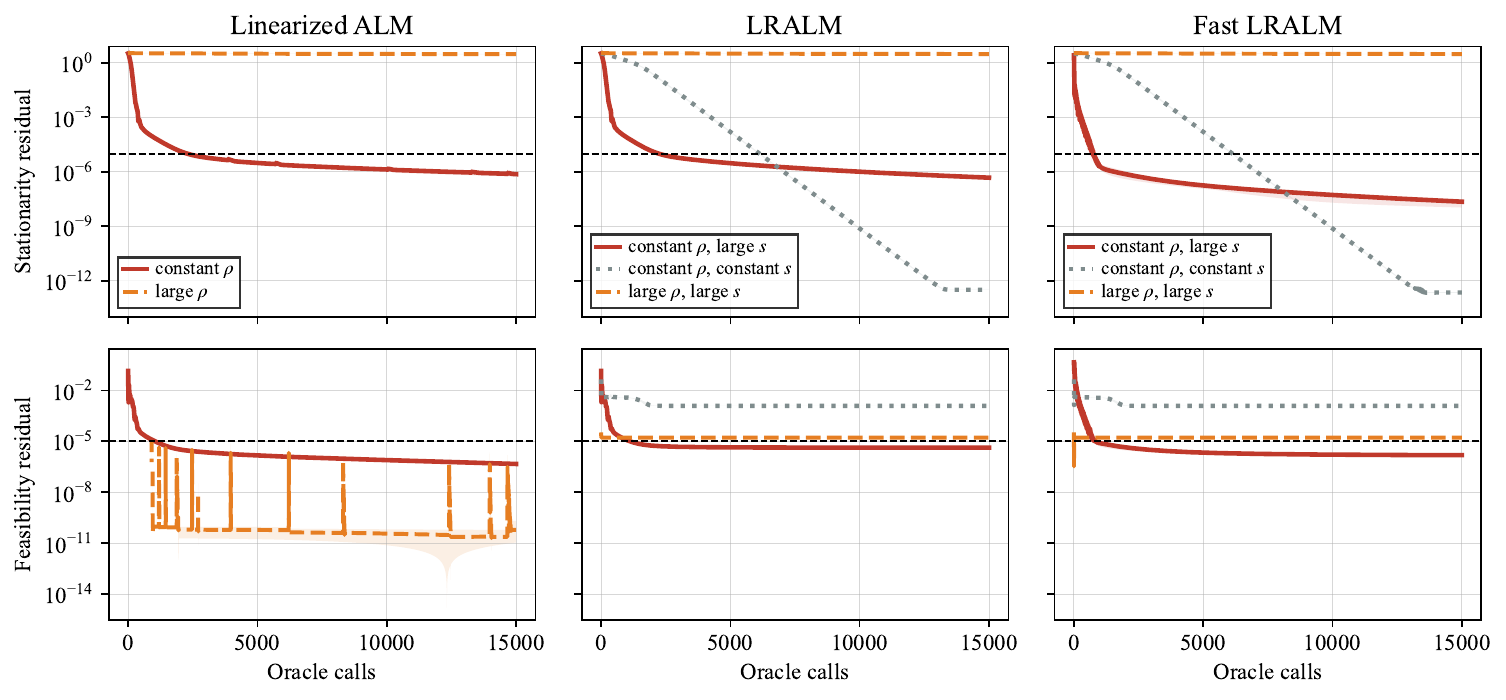}
\caption{Parameter ablation within each algorithm on Problem~\eqref{eq:exp-illcond} with \(\sigma_{\min}=10^{-4}\). The dashed line marks \(\epsilon=10^{-5}\).}
\label{fig:ablation-params-1e4}
\end{figure}

Figure~\ref{fig:ablation-algos} compares the three algorithms under constant \(\rho\), using dual regularization \(s=\Theta(\epsilon^{-1})\) for LRALM and Fast LRALM. While unaccelerated LRALM and constant-\(\rho\) linearized ALM require similar numbers of oracle calls to reach \(\epsilon\), Fast LRALM reaches \(\epsilon\) first on both instances and attains the lowest final stationarity; in the more ill-conditioned case \(\sigma_{\min}=10^{-4}\), linearized ALM reaches a lower final feasibility residual after crossing \(\epsilon\). Figures~\ref{fig:ablation-params-1e3} and~\ref{fig:ablation-params-1e4} isolate the roles of \(\rho\) and \(s\) within each algorithm. In these panels, scaling \(\rho\) as \(\Theta(\epsilon^{-1})\) stalls all three methods by shrinking the primal stepsize, whereas holding \(s\) constant leaves LRALM and Fast LRALM limited by the dual regularization bias \(\|\blmbd\|/s\), which floors the feasibility residual. These ablations show that keeping \(\rho\) constant preserves the primal stepsize, scaling \(s=\Theta(\epsilon^{-1})\) eliminates the feasibility bias, and dual Nesterov extrapolation supplies the speedup over unaccelerated methods.

\FloatBarrier

\section{Conclusion}
\label{sec:conclusion}
This paper investigates single-loop linearized ALMs with a constant penalty for
the nonsmooth nonconvex composite Problem~\eqref{prob:main}. For the LRALM
family \(\calC_\theta\), matching upper and lower bounds establish a
complexity hierarchy of \(\Theta(\epsilon^{-4})\) under
accuracy-independent parameters and \(\Theta(\epsilon^{-3})\) under
\(\Theta(\epsilon^{-1})\) dual regularization. Beyond \(\calC_\theta\), Fast
LRALM adds a Nesterov extrapolation step to the multiplier update and
attains \(\mathcal O(\epsilon^{-5/2})\) iteration complexity with a
constant penalty. To the best of our knowledge, it is the first
ALM-based single-loop method with a constant penalty to attain this rate for
Problem~\eqref{prob:main}, matching both the best-known multi-loop
complexity~\cite{kong2021accelerated} and the optimal game-stationarity
complexity for smooth NC-C minimax
problems~\cite{pan2026optimal,zhang2026matching}.
These results provide partial progress toward the open question of
Alacaoglu and Wright~\cite{Alacaoglu24a} in this composite setting.
Two natural directions remain, establishing whether
\(\mathcal O(\epsilon^{-5/2})\) is optimal for the nonconvex composite class
and extending constant-penalty ALM to nonlinearly constrained problems
as posed in~\cite{Alacaoglu24a}.

\paragraph{Acknowledgment.}
Generative artificial intelligence tools were used to assist in the
exploratory development, checking, and exposition of mathematical proofs, as
well as the implementation of numerical experiments. The authors independently
verified, critically assessed, and revised all arguments, and retain full
responsibility for the mathematical content and the final manuscript.

\setcounter{topnumber}{1}\setcounter{totalnumber}{1}

\end{document}